\documentclass[10pt]{amsart}
\usepackage{amsmath}
\usepackage{amssymb}
\usepackage{doublespace}
\usepackage[final]{graphicx}
\usepackage{float}
\usepackage{subfigure}
\usepackage{warmread}
\usepackage[authoryear]{natbib}

\usepackage[all,frame,ps,tips,line,poly,arrow,curve,import]{xy}
\SelectTips{cm}{10} \UseTips

\newcommand{\scaledfig}[2]{\scalebox{#1}{\includegraphics{#2}}}

\def\RR{\mathbb{R}}

\def\d{\mathbf{d}}

\newcommand{\bfi}{\bfseries\itshape}

\DeclareMathOperator{\SO}{SO}
\DeclareMathOperator{\SE}{SE}

\restylefloat{figure}

\newtheorem{theorem}{Theorem}[section]
\newtheorem{corollary}[theorem]{Corollary}

\newtheorem{lemma}[theorem]{Lemma}
\newtheorem{proposition}[theorem]{Proposition}

\newtheorem{remark}{Remark}[section]
\newtheorem{example}{Example}[section]

\newtheorem{definition}{Definition}[section]
\numberwithin{equation}{section}

\renewcommand{\paragraph}[1]{\vspace*{0.1in}\noindent\textbf{#1}}

\title{A Discrete Theory of Connections on Principal Bundles}
\author{Melvin Leok}
\address{Department of Mathematics, University of Michigan, Ann Arbor, MI 48109.}
\email{mleok@umich.edu}
\author{Jerrold E. Marsden}
\address{107-81, Control and Dynamical Systems, Caltech, Pasadena, CA 91125.}
\email{marsden@cds.caltech.edu}
\author{Alan D. Weinstein}
\address{Department of Mathematics, University of California, Berkeley, CA 94720.}
\email{alanw@math.berkeley.edu}

\allowdisplaybreaks

\begin{document}

\begin{abstract}
Connections on principal bundles play a fundamental role in expressing the equations of motion for mechanical systems with symmetry in an intrinsic fashion. A discrete theory of connections on principal bundles is constructed by introducing the discrete analogue of the Atiyah sequence, with a connection corresponding to the choice of a splitting of the short exact sequence. Equivalent representations of a discrete connection are considered, and an extension of the pair groupoid composition, that takes into account the principal bundle structure, is introduced. Computational issues, such as the order of approximation, are also addressed. Discrete connections provide an intrinsic method for introducing coordinates on the reduced space for discrete mechanics, and provide the necessary discrete geometry to introduce more general discrete symmetry reduction. In addition, discrete analogues of the Levi-Civita connection, and its curvature, are introduced by using the machinery of discrete exterior calculus, and discrete connections.
\end{abstract}

\maketitle

\setcounter{tocdepth}{1} \tableofcontents

\section{Introduction}

One of the major goals of geometric mechanics is the study of symmetry, and its consequences. An important tool in this regard is the non-singular reduction of mechanical systems under the action of free and proper symmetries, which is naturally formulated in the setting of principal bundles.

The reduction procedure results in the decomposition of the equations of motion into terms involving the shape and group variables, and the coupling between these are represented in terms of a connection on the principal bundle.

Connections and their associated curvature play an important role in the phenomena of geometric phases. A discussion of the history of geometric phases can be found in~\cite{Berry1990}. \cite{ShWi1989} is a collection of papers on the theory and application of geometric phases to physics.
In the rest of this section, we will survey some of the applications of geometric phases and connections to geometric mechanics and control, some of which were drawn from \cite{Marsden1994, Marsden1997, MaRa1999}.

The simulation of these phenomena requires the construction of a discrete notion of connections on principal bundles that is compatible with the approach of discrete variational mechanics, and it towards this end that this chapter is dedicated.

\begin{figure}[htbp]
\begin{center}
\includegraphics[scale=1]{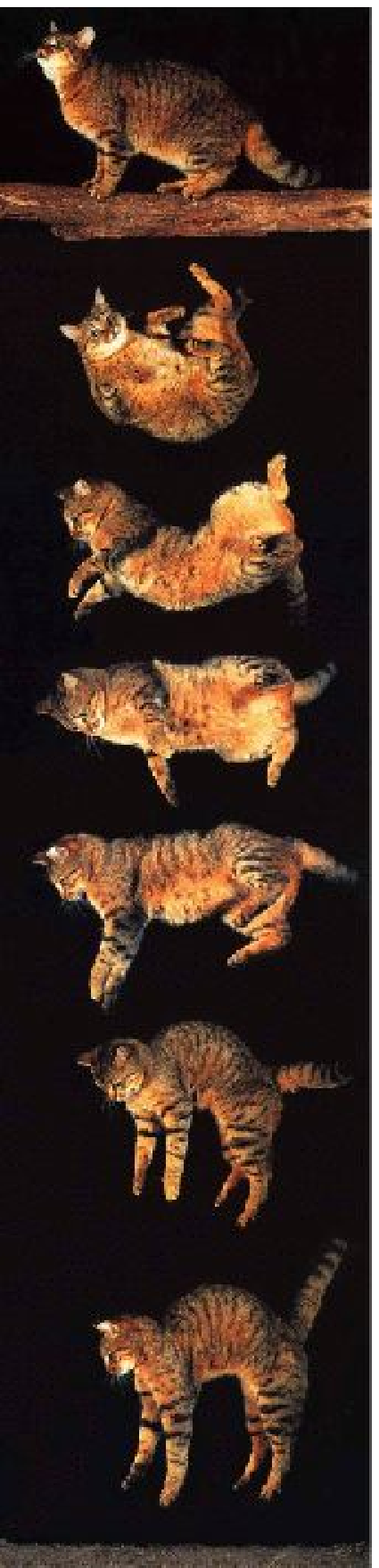}\\
\copyright\, Gerard Lacz/Animals Animals
\end{center}
\caption{\label{dcpb:fig:falling_cat}Reorientation of a falling cat at zero angular momentum.}
\end{figure}

\paragraph{Falling Cat.} Geometric phases arise in nature, and perhaps the most striking example of this is the falling cat, which is able to reorient itself by $180^\circ$, while remaining at \textit{zero} angular momentum, as show in Figure~\ref{dcpb:fig:falling_cat}.

The key to reconciling this with the constancy of the angular momentum is that angular momentum depends on the moment of inertia, which in turn depends on the shape of the cat. When the cat changes it shape by curling up and twisting, its moment of inertia changes, which is in turn compensated by its overall orientation changing to maintain the zero angular momentum condition. The zero angular momentum condition induces a connection on the principal bundle, and the curvature of this connection is what allows the cat to reorient itself.

A similar experiment can be tried on Earth, as described on page
10 of \cite{Vedral2003}. This involves standing on a swivel chair,
lifting your arms, and rotating them over your head, which will
result in the chair swivelling around slowly.

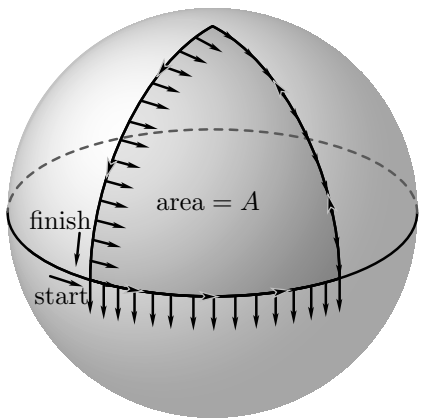
\begin{figure}[htbp]
\WARMprocessMoEPS{sphere_new}{eps}{bb}
\renewcommand{\xyWARMinclude}[1]{\scaledfig{.95}{#1}}
\begin{center}
\leavevmode
\begin{xy}
\xyMarkedImport{}
\xyMarkedTextPoints{1-2}
\xyMarkedMathPoints{3}
\end{xy}
\end{center}
\caption{\label{dcpb:fig:sphere_vector}A parallel transport of a vector around a spherical triangle produces a phase shift.}
\renewcommand{\xyWARMinclude}[1]{\scaledfig{1}{#1}}
\end{figure}

\paragraph{Holonomy.}
The sense in which curvature is related to geometric phases is most clearly illustrated by considering the parallel transport of a vector around a curve on the sphere, as shown in Figure~\ref{dcpb:fig:sphere_vector}.

Think of the point on the sphere as representing the shape of the cat, and the vector as representing its orientation. The fact that the vector experiences a phase shift when parallel transported around the sphere is an example of {\bfi holonomy}\index{holonomy}. In general, holonomy refers to a situation in geometry wherein an orthonormal frame that is parallel transported around a closed loop, back to its original position, is rotated with respect to its original orientation.

Curvature of a space is critically related to the presence of holonomy. Indeed, curvature should be thought of as being an infinitesimal version of holonomy, and this interpretation will resurface when considering the discrete analogue of curvature in the context of a discrete exterior calculus.

\begin{figure}[htbp]
\WARMprocessMoEPS{geometric_phase_foucault_pendulum_new}{eps}{bb}
\renewcommand{\xyWARMinclude}[1]{\scaledfig{.95}{#1}}
\begin{center}
\leavevmode
\begin{xy}
\xyMarkedImport{}
\xyMarkedTextPoints{1-2}
\end{xy}
\end{center}
\caption{\label{dcpb:fig:foucault_pendulum}Geometric phase of the Foucault pendulum.}
\renewcommand{\xyWARMinclude}[1]{\scaledfig{1}{#1}}
\end{figure}
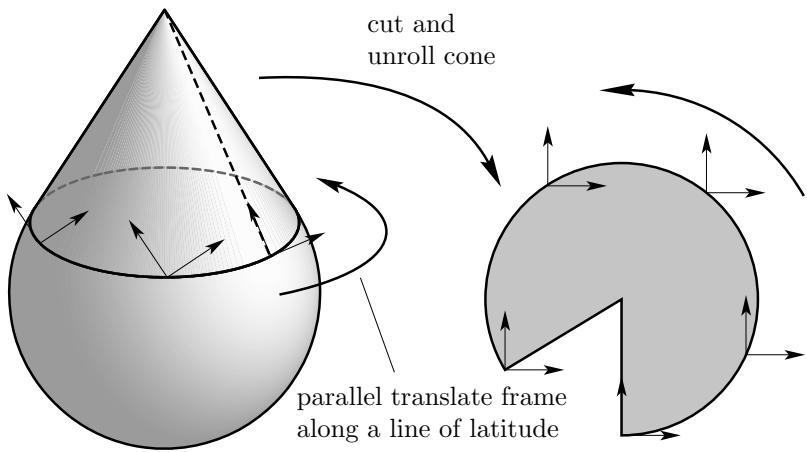

\paragraph{Foucault Pendulum.} Another example relating geometric phases and holonomy is that of the Foucault pendulum. As the Earth rotates about the Sun, the Foucault pendulum exhibits a phase shift of $\Delta \theta = 2\pi\cos\alpha$ (where $\alpha$ is the co-latitude). This phase shift is geometric in nature, and is a consequence of holonomy.  If one parallel transports an orthonormal frame around the line of constant latitude, it exhibits a phase shift that is identical to that of the Foucault pendulum, as illustrated in Figure~\ref{dcpb:fig:foucault_pendulum}.

\paragraph{True Polar Wander.} A particular striking example of the consequences of geometric phases and the conservation of angular momentum is the phenomena of true polar wander, that was studied by   \cite{GoTo1969}, and more recently by \cite{Le1998}. It is thought that some 500 to 600 million years ago, during the Vendian--Cambrian transition, the Earth, over a 15-million-year period, experienced an inertial interchange true polar wander event. This occurred when the intermediate and maximum moments of inertia crossed due to the redistribution of mass anomalies, associated with continental drift and mantle convection, thereby causing a catastrophic shift in the axis of rotation.

\begin{figure}[htbp]
\begin{center}
\includegraphics[scale=0.95]{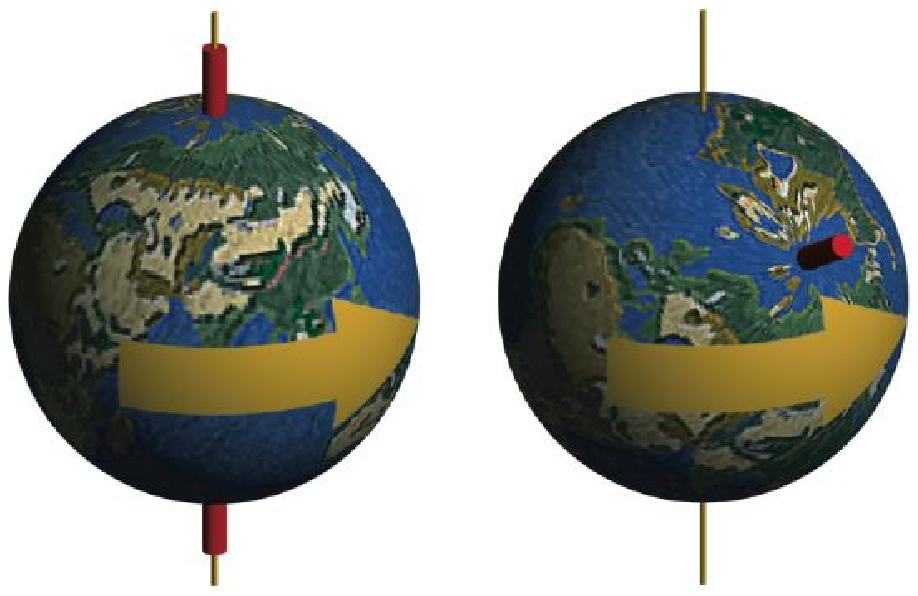}
\end{center}
\caption{\label{dcpb:fig:tpw}True Polar Wander. Red axis corresponds to the original rotational axis, and the gold axis corresponds to the instantaneous rotational axis.}
\end{figure}

This phenomena is illustrated in Figure~\ref{dcpb:fig:tpw}, wherein the places corresponding to the North and South poles of the Earth migrate towards the equator as the axis of rotation changes.

\paragraph{Geometric Control Theory.}
Geometric phases also have interesting applications and
consequences in geometric control theory, and allow, for example,
astronauts in free space to reorient themselves by changing their
shape. By holding one of their legs straight, swivelling at the
hip, and moving their foot in a circle, they are able to change
their orientation. Since the reorientation only occurs as the
shape is being changed, this allows the reorientation to be done
with extremely high precision. Such ideas have been applied to the
control of robots and spacecrafts; see, for example,
\cite{WaSa1993}. The role of connections in geometric control is
also addressed in-depth in~\cite{Marsden1994, Marsden1997}.

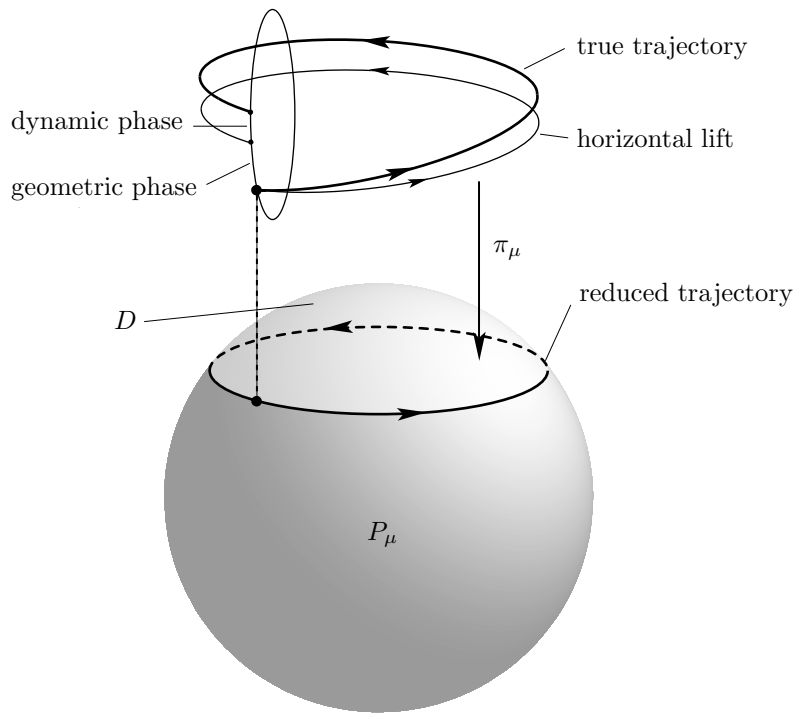
\begin{figure}[htbp]
\WARMprocessMoEPS{geometry_of_rigid_body_phase_new}{eps}{bb}
\begin{center}
\leavevmode
\begin{xy}
\xyMarkedImport{}
\xyMarkedTextPoints{1-5}
\xyMarkedMathPoints{6-8}
\end{xy}
\end{center}
\caption{\label{dcpb:fig:rigid_body_phase}Geometry of rigid-body phase.}
\end{figure}

One of the theoretical underpinnings of the application of geometric phases to geometric control was developed in \cite{Montgomery1991} and \cite{MaMoRa1990}, in the form of the {\bfi rigid-body phase formula}\index{geometric phase!rigid-body},
\[\Delta\theta = \frac{1}{\|\mu\|}\left\{\int_D\omega_\mu + 2 H_\mu T \right\} = -\Lambda +\frac{2 H_\mu T}{\|\mu\|},\]
the geometry of which is illustrated in Figure~\ref{dcpb:fig:rigid_body_phase}.

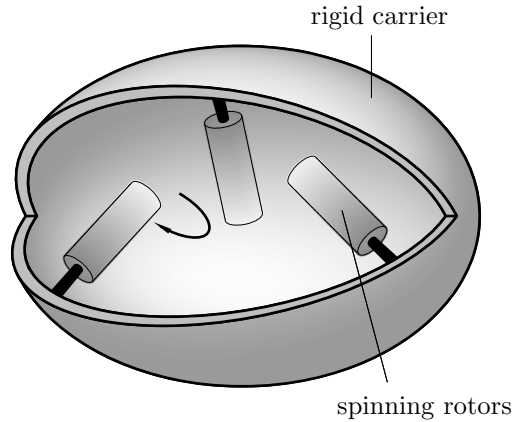
\begin{figure}[htbp]
\WARMprocessMoEPS{rigid_body_with_internal_rotors_new}{eps}{bb}
\begin{center}
\leavevmode
\begin{xy}
\xyMarkedImport{}
\xyMarkedTextPoints{1-2}
\end{xy}
\end{center}
\caption{\label{dcpb:fig:rigid_body_rotors}Rigid body with internal rotors.}
\end{figure}

An example that has been studied extensively is that of the satellite with internal rotors, with a configuration space given by $Q=\SE(3)\times S^1\times S^1\times S^1$, and illustrated in Figure~\ref{dcpb:fig:rigid_body_rotors}.

The generalization of the rigid-body phase formula in the presence of feedback control is particularly useful in the study and design of attitude control algorithms.

\section{General Theory of Bundles}

Before considering the discrete analogue of connections on principal bundles, we will review some basic material on the general theory of bundles, fiber bundles, and principal fiber bundles. A more in-depth discussion of fiber bundles can be found in \cite{Steenrod1951} and \cite{KoNo1963}.

A bundle $\mathcal{Q}$ consists of a triple $(Q,S,\pi)$, where $Q$ and $S$ are topological spaces, respectively referred to as the {\bfi bundle space}\index{bundle!bundle space} and the {\bfi base space}\index{bundle!base space}, and $\pi:Q\rightarrow S$ is a continuous map called the {\bfi projection}\index{bundle!projection}. We may assume, without loss of generality, that $\pi$ is surjective, by considering the bundle over the image $\pi(Q)\subset S$.

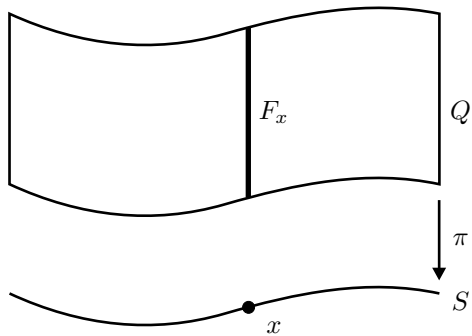
\begin{figure}[H]
\WARMprocessMoEPS{bundle}{eps}{bb}
\begin{center}
\leavevmode
\begin{xy}
\xyMarkedImport{}
\xyMarkedMathPoints{1-5}
\end{xy}
\end{center}
\caption{\label{dcpb:fig:bundle}Geometry of a fiber bundle.}
\end{figure}

The {\bfi fiber over the point $x\in S$}\index{bundle!fiber}, denoted $F_x$, is given by, $F_x=\pi^{-1}(x)$. In most situations of practical interest, the fiber at every point is homeomorphic to a common space $F$, in which case, $F$ is the {\bfi fiber} of the bundle, and the bundle is a {\bfi fiber bundle}\index{bundle!fiber bundle}. The geometry of a fiber bundle is illustrated in Figure~\ref{dcpb:fig:bundle}.

A bundle $(Q,S,\pi)$ is a {\bfi $G$-bundle}\index{bundle!$G$-bundle} if $G$ acts on $Q$ by left translation, and it is isomorphic to $(Q,Q/G,\pi_{Q/G})$, where $Q/G$ is the orbit space of the $G$ action on $Q$, and $\pi_{Q/G}$ is the natural projection.

If $G$ acts freely on $Q$, then $(Q,S,\pi)$ is called a {\bfi principal $G$-bundle}\index{bundle!principal bundle}, or {\bfi principal bundle}, and $G$ is its {\bfi structure group}\index{bundle!structure group}.  $G$ acting freely on $Q$ implies that each orbit is homeomorphic to $G$, and therefore, $\mathcal Q$ is a fiber bundle with fiber $G$.

To make the setting for the rest of this chapter more precise, we will adopt the following definition of a principal bundle,

\begin{definition}
A \textbf{principal bundle} is a  manifold $Q$ with a free
\textit{left} action, $\rho:G\times Q\rightarrow Q$, of a Lie group $G$,
such that the natural projection, $\pi:Q\rightarrow Q/G$, is a
submersion. The base space $Q/G$ is often referred to as the \textbf{shape space} $S$, which is a terminology originating from reduction theory.
\end{definition}

We will now consider a few standard techniques for combining bundles together to form new bundles. These methods include the fiber product, Whitney sum, and the associated bundle construction.

\paragraph{Fiber Product.} Given two bundles with the same base space, we can construct a new bundle, referred to as the {\bfi fiber product}\index{fiber product}, which has the same base space, and a fiber which is the direct product of the fibers of the original two bundles. More formally, we have,

\begin{definition}
Given two bundles $\pi_i:Q_i\rightarrow S,\, i=1,2$, the \textbf{fiber product} is the bundle,
\[\pi_1\times_S \pi_2: Q_1\times_S Q_2\rightarrow S,\]
where $Q_1\times_S Q_2$ is the set of all elements $(q_1,q_2)\in Q_1\times Q_2$ such that $\pi_1(q_1)=\pi_2(q_2)$, and the projection $\pi_1\times_S \pi_2$ is naturally defined by $\pi_1\times_S\pi_2(q_1,q_2)=\pi_1(q_1)=\pi_2(q_2)$. The fiber is given by $(\pi_1\times_Q\pi_2)^{-1}(x)=\pi_1^{-1}(x)\times \pi_2^{-1}(x)$.
\end{definition}

\paragraph{Whitney Sum.} The {\bfi Whitney sum}\index{Whitney sum} combines two vector bundles using the fiber product construction.

\begin{definition}
Given two vector bundles $\tau_i:V_i\rightarrow Q,\, i=1,2$, with the same base, their \textbf{Whitney sum} is their fiber product, and it is a vector bundle over $Q$, and is denoted $V_1\oplus V_2$. This bundle is obtained by taking the fiberwise direct sum of the fibers of $V_1$ and $V_2$.
\end{definition}

\paragraph{Associated Bundle.}
Given a principal bundle, $\pi:Q\rightarrow Q/G$, and a left action, $\rho:G\times M\rightarrow M$, of the Lie group $G$ on a manifold $M$, we can construct the associated bundle.

\begin{definition}\label{dcpb:def:associated_bundle}
An \textbf{associated bundle} $\tilde{M}$ with standard fiber $M$ is,
\[\tilde{M}=Q\times_G M=(Q\times M)/G,\]
where the action of $G$ on $Q\times M$ is given by $g(q,m)=(gq,gm)$. The class (or orbit) of $(q,m)$ is denoted $[q,m]_G$ or simply $[q,m]$. The projection $\pi_M:Q\times_G M\rightarrow Q/G$ is given by,
\[\pi_M:([q,m]_G)=\pi(q)\, ,\]
and it is easy to check that it is well-defined and is a surjective submersion.
\end{definition}

\section{Connections and Bundles}
Before formally introducing the precise definition of a connection, we will attempt to develop some intuition and motivation for the concept. As alluded to in the introduction to this chapter, a connection describes the curvature of a space. In the classical Riemannian setting used by Einstein in his theory of general relativity, the curvature of the space is constructed out of the connection, in terms of the Christoffel symbols that encode the connection in coordinates.

\begin{figure}[H]
\WARMprocessMoEPS{connection_new}{eps}{bb}
\begin{center}
\leavevmode
\begin{xy}
\xyMarkedImport{}
\xyMarkedTextPoints{1-6}
\end{xy}
\end{center}
\caption{\label{dcpb:fig:connections}Geometric phase and connections.}
\end{figure}
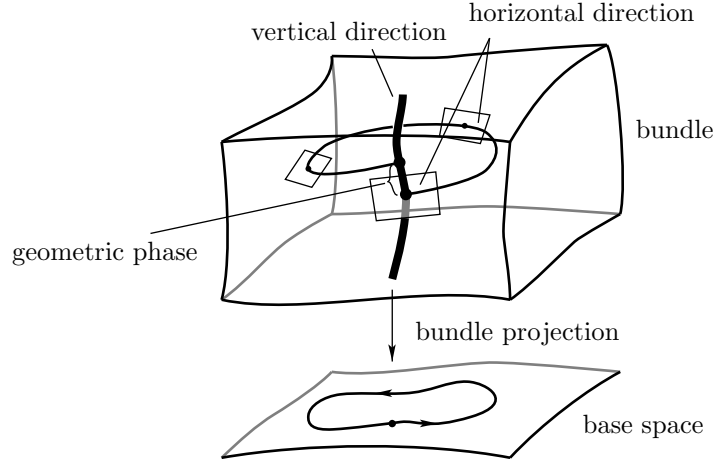

In the context of principal bundles, the connection provides a means of decomposing the tangent space to the bundle into complementary spaces, as show in Figure~\ref{dcpb:fig:connections}. Directions in the bundle that project to zero on the base space are called {\bfi vertical directions}\index{vertical!directions}, and a {\bfi connection}\index{connection} specifies a set of directions, called {\bfi horizontal directions}\index{horizontal!directions}, at each point, which complements the space of vertical directions.

In the rest of this section, we will formally define connections on principal bundles, and in the next section, discrete connections will be introduced in a parallel fashion.

\paragraph{Short Exact Sequence.}
This decomposition of the tangent space $TQ$ into horizontal and vertical subspaces yields the following short exact sequence of vector bundles over $Q$,
\[\xymatrix{
0 \ar[r] & VQ \ar[r] & TQ \ar[r]^(0.45){\pi_\ast} & \pi^* TS \ar[r] & 0\, ,
}\]
where $VQ$ is the vertical subspace of $TQ$, and $\pi^* TS$ is the pull-back of $TS$ by the projection $\pi:Q\rightarrow S$.

\paragraph{Atiyah Sequence.}
When the short exact sequence above is quotiented modulo $G$, we obtain an exact sequence of vector bundles over $S$,
\[\xymatrix{
0 \ar[r] & \tilde{\mathfrak{g}} \ar[r]^(0.35){i} & TQ/G \ar[r]^(0.6){\pi_\ast} & TS \ar[r] & 0\, ,
}\]
which is called the {\bfi Atiyah sequence}\index{Atiyah sequence} (see, for example~\cite{At1957, AlMo1985,Mackenzie1995}). Here, $\tilde{\mathfrak{g}}$ is the {\bfi adjoint bundle}\index{bundle!adjoint}, which is a special case of an associated bundle (see Definition~\ref{dcpb:def:associated_bundle}). In particular,
\[\tilde{\mathfrak{g}}=Q\times_G \mathfrak{g}=(Q\times \mathfrak{g})/G\, ,\]
where the action of $G$ on
$Q\times\mathfrak{g}$ is given by $g(q,\xi)=(gq,Ad_g\xi)$, and
$\pi_{\mathfrak{g}}:\tilde{\mathfrak{g}}\rightarrow S$ is given
by $\pi_{\mathfrak{g}}([q,\xi]_G)=\pi(q)$.

The maps in the Atiyah sequence, $i:(Q\times \mathfrak{g})/G\rightarrow TQ/G$ and $\pi_*:TQ/G\rightarrow TS$, are given by
\[i([q,\xi]_G)=[\xi_Q(q)]_G,\]
and
\[\pi_*([v_q]_G)=T\pi(v_q).\]

\paragraph{Connection $1$-form.}
Given a connection on a principal fiber bundle $\pi:Q\rightarrow Q/G$, we can represent this as a Lie algebra-valued {\bfi connection $1$-form}\index{connection!$1$-form}, $\mathcal{A}:TQ\rightarrow\mathfrak{g}$, constructed as follows (see, for example, \cite{KoNo1963}). Given an element of the Lie algebra $\xi\in\mathfrak{g}$, the infinitesimal generator map $\xi\mapsto \xi_Q$ yields a linear isomorphism between $\mathfrak{g}$ and $V_q Q$ for each $q\in Q$. For each $v_q\in T_q Q$, we define $\mathcal{A}(v_q)$ to be the unique $\xi\in\mathfrak{g}$ such that $\xi_Q$ is equal to the vertical component of $v_q$.

\begin{proposition}
The connection $1$-form, $\mathcal{A}:TQ\rightarrow\mathfrak{g}$, of a connection satisfies the following conditions.
\begin{enumerate}
\item The $1$-form is $G$-equivariant, that is,
\[\mathcal{A}\circ TL_g = \operatorname{Ad}_g \circ \mathcal{A}\, ,\]
for every $g\in G$, where $\operatorname{Ad}$ denotes the adjoint representation of $G$ in $\mathfrak{g}$.
\item The $1$-form induces a splitting of the Atiyah sequence, that is,
\[\mathcal{A}(\xi_Q)=\xi\, ,\]
for every $\xi\in\mathfrak{g}$.
\end{enumerate}
Conversely, given a $\mathfrak{g}$-valued $1$-form $\mathcal{A}$ on $Q$ satisfying conditions 1 and 2, there is a unique connection in $Q$ whose connection $1$-form is $\mathcal{A}$.
\end{proposition}
\begin{proof}
See page 64 of~\cite{KoNo1963}.
\end{proof}

\paragraph{Horizontal Lift.}
The {\bfi horizontal lift}\index{horizontal!lift} of a vector field $X\in\mathfrak{X}(S)$ is the unique vector field $X^h\in\mathfrak{X}(Q)$ which is horizontal and which projects onto $X$, that is, $T\pi_q(X^h_q)=X_{\pi(q)}$ for all $q\in Q$. The horizontal lift is in one-to-one correspondence with the choice of a connection on $Q$, as the following proposition states.

\begin{proposition}
Given a connection in $Q$, and a vector field $X\in\mathfrak{X}(S)$, there is a unique horizontal lift $X^h$ of $X$. The lift $X^h$ is left-invariant under the action of $G$. Conversely, every horizontal vector field $X^h$ on $Q$ that is left-invariant by $G$ is the lift of a vector field $X\in\mathfrak{X}(S)$.
\end{proposition}
\begin{proof}
See page 65 of~\cite{KoNo1963}.
\end{proof}

\paragraph{Connection as a Splitting of the Atiyah Sequence.}
Consider the continuous Atiyah sequence,
\[\xymatrix{
0\ar[r] & \tilde{\mathfrak{g}} \ar@<0.5ex>[r]^(0.44){i}
\ar@<-0.5ex>@{<--}[r]_(0.42){(\pi_1,\mathcal{A})}& TQ/G
\ar@<0.5ex>[r]^(0.59){\pi_*} \ar@<-0.5ex>@{<--}[r]_(0.62){X^h}
& TS \ar[r] & 0}\]
We see that the connection $1$-form, $\mathcal{A}:TQ\rightarrow\mathfrak{g}$, induces a splitting of the continuous Atiyah sequence, since
\[ (\pi_1,\mathcal{A})\circ i ([q,\xi]_G)=(\pi_1,\mathcal{A})([\xi_Q(q)]_g)=[q,\mathcal{A}(\xi_Q(q))]_G = [q,\xi]_G,\qquad\text{for all }q\in Q, \xi\in\mathfrak{g},\]
which is to say that $(\pi_1,\mathcal{A})\circ i = 1_{\tilde{\mathfrak{g}}}$. Conversely, given a splitting of the continuous Atiyah sequence, we can extend the map, by equivariance, to yield a connection $1$-form.

The horizontal lift also induces a splitting on the continuous Atiyah sequence, since, by definition, the horizontal lift of a vector field $X\in\mathfrak{X}(S)$ projects onto $X$, which is to say that $\pi_\ast\circ X^h=1_{TS}$. The horizontal lift and the connection are related by the fact that
\[ 1_{TQ/G} = i\circ(\pi_1,\mathcal{A}) + X^h \circ\pi_\ast,\]
which is a simple consequence of the fact that the two splittings are part of the following commutative diagram,
\[\xymatrix{
0\ar[r] & \tilde{\mathfrak{g}}
\ar@{=}[d]_{1_{\tilde{\mathfrak{g}}}} \ar@<0.5ex>[r]^(0.44){i}
\ar@<-0.5ex>@{<--}[r]_(0.42){(\pi_1,\mathcal{A})}& TQ/G
\ar@<0.5ex>[r]^(0.59){\pi_*} \ar@<-0.5ex>@{<--}[r]_(0.62){X^h}
\ar[d]^{\alpha_{\mathcal{A}}}& TS \ar@{=}[d]^{1_{TS}}
\ar[r] & 0\\
0\ar[r] & \tilde{\mathfrak{g}} \ar@<0.5ex>[r]^(0.4){i_1}
\ar@<-0.5ex>@{<--}[r]_(0.4){\pi_1} & \tilde{\mathfrak{g}}\oplus TS
\ar@<0.5ex>[r]^(0.6){\pi_2} \ar@<-0.5ex>@{<--}[r]_(0.6){i_2}& TS
\ar[r] & 0 }\]
where $\alpha_{\mathcal{A}}$ is an isomorphism. The isomorphism is given in the following lemma.

\begin{lemma}
The map $\alpha_{\mathcal{A}}:TQ/G\rightarrow \tilde{\mathfrak{g}}\oplus TS$ defined by
\[ \alpha_\mathcal{A}([q,\dot q]_G)=[q,\mathcal{A}(q,\dot q)]_G\oplus T\pi(q,\dot q), \]
is a well-defined vector bundle isomorphism. The inverse of $\alpha_{\mathcal{A}}$ is given by
\[\alpha^{-1}_{\mathcal{A}}([q,\xi]_G\oplus (x,\dot x))=[(x,\dot x)^h_q+\xi q]_G.\]
\end{lemma}
\begin{proof}
See page 15 of~\cite{CeMaRa2001}.
\end{proof}

This lemma, and its higher-order generalization, that identifies $T^{(2)} Q/G$ with $T^{(2)}S\times_{S}2\tilde{\mathfrak{g}}$, is critical in allowing us to construct the Lagrange--Poincar\'e operator, which is an intrinsic method of expressing the reduced equations arising from Lagrangian reduction.

In the next section, we will develop the theory of discrete connections on principal bundles in a parallel fashion to the way we introduced continuous connections.

\section{Discrete Connections}

Discrete variational mechanics is based on the idea of approximating the tangent bundle $TQ$ of Lagrangian mechanics with the pair groupoid $Q\times Q$. As such, the purpose of a discrete connection is to decompose the subset of $Q\times Q$ that projects to a neighborhood of the diagonal of $S\times S$ into horizontal and vertical spaces.

The reason why we emphasize that the construction is only valid for the subset of $Q\times Q$ that projects to a neighborhood of the diagonal of $S\times S$ is that there are topological obstructions to globalizing the construction to all of $Q\times Q$ except in the case that $Q$ is a trivial bundle.

One of the challenges of dealing with the discrete space modelled
by the pair groupoid $Q\times Q$ is that it is not a linear space,
in contrast to $TQ$. As we shall see, the standard pair groupoid
composition is not sufficient to make sense of the notion of an
element $(q_0,q_1)\in Q\times Q$ being the composition of a
horizontal and a vertical element. We will propose a natural
notion of composing an element with a vertical element that
makes sense of the horizontal and vertical decomposition.

In the subsequent sections, we will use the discrete connection to extend the pair groupoid composition even further, and explore its applications to the notion of curvature in discrete geometry.

\begin{figure}[H]
\begin{center}
\includegraphics{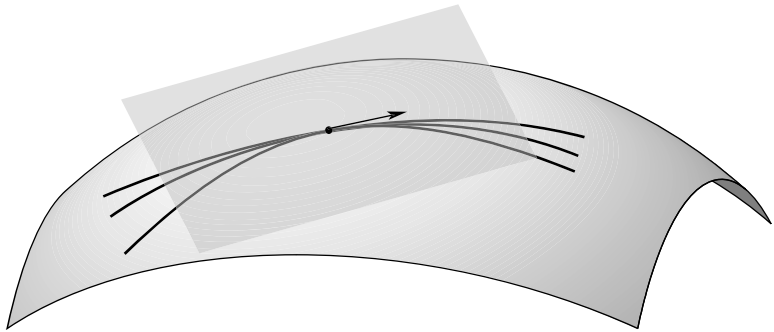}
\end{center}
\caption{\label{dcpb:fig:curvesvectors}Intrinsic representation of the tangent bundle.}
\end{figure}

\paragraph{Intrinsic Representation of the Tangent Bundle.}
The intuition underlying our construction of discrete horizontal and vertical spaces is best developed by considering the intrinsic representation of the tangent bundle. This representation is obtained by identifying a tangent vector at a point on the manifold with the equivalence class of curves on the manifold going through the point, such that the tangent to the curve at the point is given by the tangent vector. This notion is illustrated in Figure~\ref{dcpb:fig:curvesvectors}.

Given a vector $v_q\in TQ$, we identify it with the family of curves $q:\mathbb{R}\rightarrow Q$, such that $q(0)=q$, and $\dot q(0)=v$. The equivalence class $[\,\cdot\,]$ identifies curves with the same basepoint, and the same velocity at the basepoint.

With this representation, it is natural to consider $(q_0,q_1)\in Q\times Q$ to be an approximation of $[q(\cdot)]=v_q\in TQ$, in the sense that,
\[
q_0=q(0),\qquad\qquad
q_1=q(h),
\]
for some fixed time step $h$, and where $q(\cdot)$ is a representative curve corresponding to $v_q$ in the intrinsic representation of the tangent bundle.

\subsection{Horizontal and Vertical Subspaces of $Q\times Q$}\label{dcpb:subsec:hor_ver}
Recall that the vertical subspace at a point $q$, denoted $V_q$, is given by
\[V_q=\{v_q\in TQ\mid \pi_*(v_q)=0\}=\{\xi_Q\mid
\xi\in\mathfrak{g}\}.\]
Notice that the vertical space is precisely that subspace of $TQ$ which maps under the lifted projection map to the embedded copy of $S$ in $TS$. We proceed in an analogous fashion to define a discrete vertical subspace at a point $q$.

The natural discrete analogue of the lifted projection map $\pi_*$ is the diagonal action of the projection map on $Q\times Q$, $(\pi,\pi):Q\times Q\rightarrow Q\times Q$, where $(q_0,q_1)\mapsto(\pi q_0, \pi q_1)$. This is because
\[ \pi_*(v_q) = \pi_*([q(\cdot)]) = [\pi(q(\cdot))].\]
In the same way that we embed $S$ into $TS$ by the map $x\mapsto [x]=0_x$, $S$ naturally embeds itself into the diagonal of $S\times S$, $x\mapsto(x,x)=e_{S\times S}$, which we recall is the identity subspace of the pair groupoid.

The alternative description of the vertical space is in terms of the embedding of $Q\times\mathfrak{g}$ into $TQ$, by $(q,\xi)\mapsto \xi_Q(q)$, using the infinitesimal generator construction,
\[ \xi_Q(q) = [ \exp(\xi t) q ].\]
In an analogous fashion, we construct a {\bfi discrete generator}\index{discrete generator} map, which is given in the following definition.

\begin{definition}
The \textbf{discrete generator} is the map $i:Q\times G\rightarrow Q\times Q$, given by
\[ i(q,g) = (q,gq),\]
which we also denote by $i_q(g) = i(q,g) = (q,gq)$.
\end{definition}

Then, we have the following definition of the {\bfi discrete vertical space}\index{vertical!space!discrete}.
\begin{definition}
The \textbf{discrete vertical space} is given by
\begin{align*} \operatorname{Ver}_q
&=\{(q,q')\in Q\times Q \mid
(\pi,\pi)(q,q')=e_{S\times S}\}\\
&=\{i_q(g)\mid g\in G\}.
\end{align*}
This is the discrete analogue of the statement
$\operatorname{Ver}_q=\{v_q\in TQ\mid \pi_*(v_q)=0\}=\{\xi_Q\mid
\xi\in\mathfrak{g}\}$.
\end{definition}

Since the pair groupoid composition is only defined on the space of composable pairs, we need to extend the composition to make sense of how the discrete horizontal space is complementary to the discrete vertical space. In particular, we define the composition of a vertical element with an arbitrary element of $Q\times Q$ as follows.

\begin{definition}
The composition of an arbitrary element $(q_0,q_1)\in Q\times Q$
with a vertical element is given by
\[ i_{q_0}(g) \cdot (q_0,q_1) =
(e,g)(q_0,q_1)=(q_0,gq_1).\]
\end{definition}

An elementary consequence of this definition is that it makes the discrete generator map a homomorphism.

\begin{lemma}\label{dcpb:lemma:discrete_generator_homomorphism}
The discrete generator, $i_q$, is a homomorphism. This is a discrete analogue of the statement in the continuous theory that $(\xi+\chi)_Q = \xi_Q +\chi_Q$.
\end{lemma}
\begin{proof}
We compute,
\begin{align*}
i_q(g)\cdot i_q(h)&=i_q(g)\cdot(q,hq)\\
&=(e,g)(q,hq)\\
&=(q,ghq)\\
&=i_q(gh).
\end{align*}
Therefore, $i_q$ is a homomorphism.
\end{proof}

If we define the $G$ action on $Q\times G$ to be $h(q,g)=(hq, hgh^{-1})$, we find that the composition of a vertical element with an arbitrary element is $G$-equivariant.

\begin{lemma}\label{dcpb:lemma:vertical_composition_equivariant}
The composition of a vertical element with an arbitrary element of $Q\times Q$ is $G$-equivariant,
\[ i_{hq_0}(hgh^{-1}) \cdot (hq_0,hq_1)  = h\cdot i_{q_0}(g)\cdot(q_0,q_1)\, . \]
\end{lemma}
\begin{proof}
Consider the following computation,
\begin{align*}
i_{hq_0}(hgh^{-1}) \cdot (hq_0,hq_1)
&= (hq_0, hgh^{-1} h q_1)\\
&= (hq_0, hgq_1)\\
&= h(q_0, gq_1)\\
&= h\cdot i_{q_0}(g)\cdot (q_0,q_1)\,.\qedhere
\end{align*}
\end{proof}

Having made sense of how to compose an arbitrary element of
$Q\times Q$ with a vertical element, we are in a position to
introduce the notion of a discrete connection.

A {\bfi discrete connection}\index{connection!discrete|see{discrete connection}} is a $G$-equivariant choice of a subset of $Q\times Q$ called the {\bfi discrete horizontal space}\index{horizontal!space!discrete}, that is complementary to the discrete vertical space. In particular, given $(q_0,q_1)\in Q\times Q$, a discrete connection decomposes this into the {\bfi horizontal component}\index{horizontal!component}, $\operatorname{hor}(q_0,q_1)$, and the {\bfi vertical component}\index{vertical!component}, $\operatorname{ver}(q_0,q_1)$, such that
\[ \operatorname{ver}(q_0,q_1)\cdot\operatorname{hor}(q_0,q_1) = (q_0,q_1)\, ,\]
in the sense of the composition of a vertical element with an arbitrary element we defined previously. Furthermore, the $G$-equivariance condition states that
\begin{align*}
\operatorname{hor}(g q_0, g q_1) &= g \cdot \operatorname{hor}(q_0,q_1) \, ,\\
\intertext{and}
\operatorname{ver}(g q_0, g q_1) &= g \cdot \operatorname{ver}(q_0,q_1) \, .
\end{align*}

\subsection{Discrete Atiyah Sequence}
Recall that we obtain a short exact sequence corresponding to the decomposition of $TQ$ into horizontal and vertical spaces. Due to the equivariant nature of the decomposition, quotienting this short exact sequence yields the Atiyah sequence. In this subsection, we will introduce the analogous discrete objects.

\paragraph{Short Exact Sequence.}
The decomposition of the pair groupoid $Q\times Q$, into discrete
horizontal and vertical spaces, yields the following short exact
sequence of bundles over $Q$.
\[\xymatrix{
0 \ar[r] & \operatorname{Ver} Q \ar[r]^(0.47){i} & Q\times Q \ar[r]^(0.40){(\pi,\pi)} & (\pi,\pi)^* S\times S \ar[r] & 0,
}\]
where $\operatorname{Ver} Q$ is the discrete vertical subspace of $Q\times Q$, and $(\pi,\pi)^* S\times S$ is the pull-back of $S\times S$ by the projection $(\pi,\pi):Q\times Q\rightarrow S\times S$.

\paragraph{Discrete Atiyah Sequence.}
When the short exact sequence above is quotiented modulo $G$, we obtain an exact sequence of bundles over $S$,
\[\xymatrix{
0\ar[r]
& \tilde{G} \ar[r]^(0.34){i}
& (Q\times Q)/G \ar[r]^(0.59){(\pi,\pi)}
& S \times S \ar[r]
& 0\, ,}\]
which we call the {\bfi discrete Atiyah sequence}\index{Atiyah sequence!discrete|see{discrete Atiyah sequence}}\index{discrete Atiyah sequence}. Here, $\tilde{G}$ is an associated bundle (see Definition~\ref{dcpb:def:associated_bundle}). In particular,
\[ \tilde{G} = Q\times_G G = (Q\times G)/G\, ,\]
where the action of $G$ on $Q\times G$ is given by $g(q,h) = (gq, ghg^{-1})$, which is the natural discrete analogue of the adjoint action of $\mathfrak{g}$ on $Q\times\mathfrak{g}$. Furthermore, $\pi_G:\tilde{G}\rightarrow S$ is given by $\pi_G([q,g]_G)=\pi(q)$.

The maps in the discrete Atiyah sequence $i:\tilde{G} \rightarrow (Q\times Q)/G$, and $(\pi,\pi):(Q\times Q)/G]\rightarrow S\times S$, are given by
\[ i([q,g]_G)=[q,gq]_G=[i_q(g)]_G\, ,\]
and
\[ (\pi,\pi)([q_0,q_1]_g)=(\pi q_0,\pi q_1)\, .\]

\subsection{Equivalent Representations of a Discrete Connection}

In addition to the discrete connection which arises from a $G$-equivariant decomposition of the pair groupoid $Q\times Q$ into a discrete horizontal and vertical space, we have equivalent representations in terms of splittings of the discrete Atiyah sequence, as well as maps on the unreduced short exact sequence.

\paragraph{Maps on the Unreduced Short Exact Sequence.} These correspond to discrete analogues of the connection $1$-form, and the horizontal lift.
\begin{itemize}
\item Discrete connection $1$-form, $\mathcal{A}_d:Q\times Q \rightarrow G$.
\item Discrete horizontal lift, $(\cdot,\cdot)^h_q:S\times Q\rightarrow Q\times Q$.
\end{itemize}

\paragraph{Maps That Yield a Splitting of the Discrete Atiyah Sequence.}
\begin{itemize}
\item $(\pi_1,\mathcal{A}_d):(Q\times Q)/G\rightarrow \tilde{G}$, which is related to the discrete connection $1$-form.
\item $(\cdot,\cdot)^h:S\times S\rightarrow (Q\times Q)/G$, which is related to the discrete horizontal lift.
\end{itemize}

\paragraph{Relating the Two Sets of Representations.}
These two sets of representations are related in the following way:
\begin{itemize}
\item The maps on the unreduced short exact sequence are equivariant, and hence drop to the discrete Atiyah sequence, where they induce splittings of the short exact sequence.
\item The maps that yield splittings of the discrete Atiyah sequence can be extended equivariantly to recover the maps on the unreduced short exact sequence.
\end{itemize}
Furthermore, standard results from homological algebra yield an equivalence between the two splittings of the discrete Atiyah sequence.

In the rest of this section, we will also discuss in detail the method of moving between the various representations of the discrete connection. The organization of the rest of the section, and the subsections in which we relate the various representations are given in the following diagram.
\[
\xymatrix@R=2.5pc@C=-3pc{
& *=<5cm,1.6cm>+[F-,]+{\txt{\S\ref{dcpb:subsec:hor_ver}\\ discrete connection\\ $\operatorname{hor}:Q\times Q\rightarrow \operatorname{Hor}_q$\\ $\operatorname{ver}:Q\times Q\rightarrow\operatorname{Ver}_q$}} &\\
*=<5cm,1.25cm>+[F-,]+{\txt{\S\ref{dcpb:subsec:discrete_connection_one_form}\\ discrete connection $1$-form\\ $\mathcal{A}_d:Q\times Q\rightarrow G$}} \ar@{<->}[ur]^{\S\ref{dcpb:subsec:discrete_connection_one_form}}& &
*=<5cm,1.25cm>+[F-,]+{\txt{\S\ref{dcpb:subsec:discrete_horizontal_lift}\\ discrete horizontal lift\\
$(\cdot,\cdot)^h_q:S\times S\rightarrow Q\times Q$}} \ar@{<->}[ul]_{\S\ref{dcpb:subsec:discrete_horizontal_lift}} \ar@{<->}[ll]_{\S\ref{dcpb:subsec:discrete_horizontal_lift}}\\
*=<5cm,1.25cm>+[F-,]+{\txt{\S\ref{dcpb:subsec:splitting_connection}\\ splitting (connection $1$-form)\\ $(\pi,\mathcal{A}_d):(Q\times Q)/G\rightarrow \tilde{G}$}}
\ar@{<->}[u]^{\S\ref{dcpb:subsec:splitting_connection}}& &
*=<5cm,1.25cm>+[F-,]+{\txt{\S\ref{dcpb:subsec:splitting_lift}\\ splitting (horizontal lift) \\
$(\cdot,\cdot)^h:S\times S\rightarrow (Q\times Q)/G$}}
\ar@{<->}[u]_{\S\ref{dcpb:subsec:splitting_lift}}
}
\]\index{discrete connection!relating the representations}

\subsection{Discrete Connection $1$-Form}\label{dcpb:subsec:discrete_connection_one_form}
Given a discrete connection on a principal fiber bundle $\pi:Q\rightarrow Q/G$, we can represent this as a Lie group-valued {\bfi discrete connection $1$-form}\index{discrete connection!$1$-form|see{discrete connection $1$-form}}, $\mathcal{A}_d:Q\times Q\rightarrow G$, which is a natural generalization of the Lie algebra-valued connection $1$-form on tangent bundles, $\mathcal{A}:TQ\rightarrow\mathfrak{g}$, to the discrete context.

\paragraph{Discrete Connection $1$-Forms from Discrete Connections.}\index{discrete connection $1$-form!from discrete connection}\index{discrete connection!to discrete connection $1$-form}
The discrete connection $1$-form is constructed as follows. Given an element of the Lie group $g\in G$, the discrete generator map $g\mapsto i_q(g)$ yields an isomorphism between $G$ and $\operatorname{Ver}_q$ for each $q\in Q$. For each $(q_0,q_1)\in Q\times Q$, we define $\mathcal{A}_d(q_0,q_1)$ to be the unique $g\in G$ such that $i_q(g)$ is equal to the vertical component of $(q_0,q_1)$. In particular, this is equivalent to the condition that the following statement holds,
\[ (q_0,q_1)=i_{q_0} (\mathcal{A}_d(q_0,q_1)) \cdot \operatorname{hor}(q_0,q_1) .\]

\begin{remark}
It follows from the above identity that the discrete horizontal space can also be expressed as
\begin{align*}
\operatorname{Hor}_{q_0}
&=\{ (q_0,q_1)\in Q\times Q \mid \operatorname{hor}(q_0,q_1)=(q_0,q_1)\}\\
&=\{ (q_0,q_1)\in Q\times Q \mid \mathcal{A}_d(q_0,q_1)=e\}\, .
\end{align*}
\end{remark}

We will now establish a few properties of the discrete connection $1$-form.

\begin{proposition}
The \textbf{discrete connection $1$-form}\index{discrete connection $1$-form!properties}, $\mathcal{A}_d:Q\times
Q\rightarrow G$, satisfies the following properties.
\begin{enumerate}
  \item The $1$-form is $G$-equivariant, that is, \[\mathcal{A}_d\circ
  L_g=I_g\circ\mathcal{A}_d,\] which is the discrete analogue of
  the $G$-equivariance of the continuous connection, $\mathcal{A}\circ
  TL_g=Ad_g\circ\mathcal{A}$.
  \item The $1$-form induces a splitting of the Discrete Atiyah
  sequence, that is, \[\mathcal{A}_d(i_q(g))=\mathcal{A}_d(q_0,gq_0)=g,\] which is the discrete analogue of
$\mathcal{A}(\xi_Q)=\xi$.
\end{enumerate}
\end{proposition}
\begin{proof}
The proof relies on the properties of a discrete connection, and the definition of the discrete connection $1$-form.
\begin{enumerate}
\item
The discrete connection $1$-form satisfies the condition
\[ (q_0,q_1) =  i_{q_0}(\mathcal{A}_d(q_0,q_1)) \cdot \operatorname{hor}(q_0,q_1)\, .\]
If we denote $\operatorname{hor}(q_0,q_1)$ by $(q_0,\bar q_1)$, we have that
\[ (q_0,q_1) = (q_0, \mathcal{A}_d(q_0,q_1) \bar q_1)\, .\]
Similarly, we have,
\begin{align*}
(gq_0, gq_1)
&= i_{gq_0}(\mathcal{A}_d(gq_0,gq_1))\cdot\operatorname{hor}(gq_0, gq_1)\\
&= i_{gq_0}(\mathcal{A}_d(gq_0,gq_1))\cdot g\cdot \operatorname{hor}(q_0,q_1)\\
&= (e,\mathcal{A}_d(gq_0,gq_1))(gq_0,g\bar q_1)\\
&= (gq_0, \mathcal{A}_d(gq_0,gq_1)g\bar q_1)\, ,
\end{align*}
where we have used the $G$-equivariance of the discrete horizontal space. By looking at the expressions for $gq_1$ and $q_1$, we conclude that
\begin{align*}
\mathcal{A}_d(gq_0,gq_1) g \bar q_1
&= gq_1\\*
&= g \mathcal{A}_d(q_0,q_1) \bar q_1\, ,\\
\mathcal{A}_d(gq_0,gq_1) g &= g \mathcal{A}_d(q_0,q_1)\, ,\\
\mathcal{A}_d(gq_0,gq_1) &= g \mathcal{A}_d(q_0,q_1) g^{-1}\, ,
\end{align*}
which is precisely the statement that $\mathcal{A}_d\circ L_g  = I_g \circ \mathcal{A}_d$, that is to say that $\mathcal{A}_d$ is $G$-equivariant.
\item
Recall that $i_q(g)$ is an element of the discrete vertical space. Since the discrete horizontal space is complementary to the discrete vertical space, it follows that $\operatorname{ver}(i_q(g))=i_q(g)$. Then, by the construction of the discrete connection $1$-form, $\mathcal{A}_d(i_q(g))$ is the unique element of $G$ such that
\[ i_q(\mathcal{A}_d(i_q(g)))=\operatorname{ver}(i_q(g))=i_q(g)\,.\]
Since $i_q$ is an isomorphism between $G$ and the discrete
vertical space, we conclude that $\mathcal{A}_d(i_q(g)) = g$, as
desired.\qedhere
\end{enumerate}
\end{proof}

The second result is equivalent to the map recovering the discrete Euler--Poincar\'{e} connection when restricted to a $G$-fiber, that is, $\mathcal{A}_d(x,g_0,x,g_1)=g_1 g_0^{-1}.$ In particular, it follows that the map is trivial when restricted to the diagonal space, that is, $\mathcal{A}_d(q,q)=e$.

The properties of a discrete connection are discrete analogues of the properties of a continuous connection in the sense that if a discrete connection has a given property, the corresponding continuous connection which is induced in the infinitesimal limit has the analogous continuous property. The precise sense in which a discrete connection induces a continuous connection will be discussed in \S\ref{s:infinitesimal_limit}.

\paragraph{Discrete Connections from Discrete Connection $1$-Forms.}\index{discrete connection!from discrete connection $1$-form}\index{discrete connection $1$-form!to discrete connection}
Having shown how to obtain a discrete connection $1$-form from a discrete connection, let us consider the converse case of obtaining a discrete connection from a discrete connection $1$-form with the properties above. We do this by constructing the discrete horizontal and vertical components as follows.

\begin{definition}
Given a discrete connection $1$-form, $\mathcal{A}_d:Q\times Q \rightarrow G$ that is $G$-equivariant and induces a splitting of the discrete Atiyah sequence, we define the \textbf{horizontal component} to be
\[ \operatorname{hor}(q_0,q_1) = i_{q_0}((\mathcal{A}_d(q_0,q_1))^{-1}) \cdot (q_0,q_1)\, .\]
The \textbf{vertical component} is given by
\[ \operatorname{ver}(q_0,q_1) = i_{q_0} (\mathcal{A}_d(q_0,q_1))\, .\]
\end{definition}

\begin{proposition}
The discrete connection we obtain from a discrete connection $1$-form has the following properties.
\begin{enumerate}
\item The discrete connection yields a horizontal and vertical decomposition of $Q\times Q$, in the sense that
\[ (q_0,q_1)= \operatorname{ver}(q_0,q_1)\cdot\operatorname{hor}(q_0,q_1)\, , \]
for all $(q_0,q_1)\in Q\times Q$.
\item The discrete connection is $G$-equivariant, in the sense that
\begin{align*}
\operatorname{hor}(g q_0, g q_1) &= g \cdot \operatorname{hor}(q_0,q_1) \, ,\\
\intertext{and}
\operatorname{ver}(g q_0, g q_1) &= g \cdot \operatorname{ver}(q_0,q_1) \, .
\end{align*}
\end{enumerate}
\end{proposition}
\begin{proof}
The proof relies on the properties of the discrete connection $1$-form, and the definitions of the discrete horizontal and vertical spaces.
\begin{enumerate}
\item Consider the following computation,
\begin{align*}
\operatorname{ver}(q_0,q_1)\cdot\operatorname{hor}(q_0,q_1)
&=i_{q_0} (\mathcal{A}_d(q_0,q_1)) \cdot i_{q_0} ((\mathcal{A}_d(q_0,q_1))^{-1}) \cdot (q_0,q_1)\\
&= i_{q_0}(\mathcal{A}_d(q_0,q_1)(\mathcal{A}_d(q_0,q_1))^{-1}) \cdot (q_0,q_1)\\
&= i_{q_0}(e)\cdot(q_0,q_1)\\
&= (q_0,q_1)\, ,
\end{align*}
where we used that $i_q$ is a homomorphism (see Lemma~\ref{dcpb:lemma:discrete_generator_homomorphism}).
\item We compute,
\begin{align*}
\operatorname{hor}(g q_0, g q_1)
&= i_{gq_0}((\mathcal{A}_d(gq_0,gq_1))^{-1}) \cdot (g q_0, g q_1)\\
&= i_{gq_0}(g(\mathcal{A}_d(q_0,q_1))^{-1}g^{-1}) \cdot (g q_0, g q_1)\\
&= (e,g(\mathcal{A}_d(q_0,q_1))^{-1}g^{-1})(g q_0, g q_1)\\
&= (g q_0, g (\mathcal{A}_d(q_0,q_1))^{-1} q_1)\\
&= g \cdot i_{q_0} ((\mathcal{A}_d(q_0,q_1))^{-1})\cdot (q_0,q_1)\\
&= g \cdot \operatorname{hor}(q_0,q_1)\, ,
\end{align*}
where we have used the fact that the composition of a vertical element with an arbitrary element is $G$-equivariant (see Lemma~\ref{dcpb:lemma:vertical_composition_equivariant}). Similarly, we compute,
\begin{align*}
\operatorname{ver}(gq_0,gq_1)
&= i_{gq_0} (\mathcal{A}_d(gq_0,gq_1))\\
&= i_{gq_0} (g\mathcal{A}_d(q_0,q_1)g^{-1})\\
&= (gq_0, g\mathcal{A}_d(q_0,q_1)g^{-1}gq_0)\\
&= g \cdot (q_0, \mathcal{A}_d(q_0,q_1)q_0)\\
&= g \cdot i_{q_0}(\mathcal{A}_d(q_0,q_1))\\
&=  g \cdot \operatorname{ver}(q_0,q_1)\, .\qedhere
\end{align*}
\end{enumerate}
\end{proof}

\paragraph{Local Representation of the Discrete Connection $1$-Form.}\index{discrete connection $1$-form!local representation}
Since the discrete connection $1$-form can be thought of as
comparing group fiber quantities at different base points, we obtain
the natural identity that
\[\mathcal{A}_d(gq_0,hq_1)=h\mathcal{A}_d(q_0,q_1)g^{-1}.\]
In a local trivialization, this corresponds to
\[\mathcal{A}_d(x_0,g_0,x_1,g_1)=g_1\mathcal{A}_d(x_0,e,x_1,e)g_0^{-1}.\]
We define
\[A(x_0,x_1)=\mathcal{A}_d(x_0,e,x_1,e),\]
which yields the local representation of the discrete connection
$1$-form.
\begin{definition}\label{dcpb:defn:local_representation}
Given a discrete connection $1$-form, $\mathcal{A}_d:Q\times Q\rightarrow G$, its \textbf{local representation} is given by
\[\mathcal{A}_d(x_0,g_0,x_1,g_1)=g_1 A(x_0,x_1) g_0^{-1}\, ,\]
where
\[A(x_0,x_1)=\mathcal{A}_d(x_0,e,x_1,e)\, .\]
\end{definition}

\begin{lemma}
The local representation of a discrete connection is $G$-equivariant.
\end{lemma}
\begin{proof}
Consider the following computation,
\begin{align*}
\mathcal{A}_d(g(x_0,g_0),g(x_1,g_1))
&= \mathcal{A}_d((x_0, gg_0),(x_1,gg_1))\\
&= gg_1 A(x_0,x_1)(gg_0)^{-1}\\
&= g(g_1 A(x_0,x_1)g_0^{-1})g^{-1}\\
&=g \mathcal{A}_d((x_0,g_0),(x_1,g_1)) g^{-1}\, ,
\end{align*}
which shows that the local representation is $G$-equivariant, as expected.
\end{proof}

Notice also that in the pure group case, where $Q=G$, this recovers the discrete Euler--Poincar\'e connection, as we would expect, since the shape space is trivial. In particular,  $x_0=x_1=e$, which implies that $A(x_0,x_1)=\mathcal{A}_d(e,e,e,e)=e$, and $\mathcal{A}_d(g_0,g_1)=g_1 g_0^{-1}$.

\begin{example}\index{discrete connection!mechanical}
As an example, we construct the natural discrete analogue of the
mechanical connection, $\mathcal{A}:TQ\rightarrow\mathfrak{g}$, by
the following procedure, which yields a discrete connection $1$-form, $\mathcal{A}_d:Q\times Q \rightarrow G$.
\begin{enumerate}
  \item Given the point $(q_0,q_1)\in Q\times Q$, we construct the
  geodesic path $q_{01}:[0,1]\rightarrow Q$ with respect to the
  kinetic energy metric, such that $q_{01}(0)=q_0$, and
  $q_{01}(1)=q_1$.
  \item Project the geodesic path to the shape space, $x_{01}(t)\equiv\pi
  q_{01}(t)$, to obtain the curve $x_{01}$ on $S$.
  \item Taking the horizontal lift of $x_{01}$ to $Q$ using the connection
  $\mathcal{A}$ yields $\tilde{q}_{01}$.
  \item There is a unique $g\in G$ such that
  $q_{01}(1)=g\cdot \tilde{q}_{01}(1)$.
  \item Define $\mathcal{A}_d(q_0,q_1)=g$.
\end{enumerate}
This discrete connection is consistent with the classical notion
of a connection in the limit that $q_1$ approaches $q_0$, in the
usual sense in which discrete mechanics on $Q\times Q$ converges
to continuous Lagrangian mechanics on $TQ$. As mentioned before, this statement is made more precise in \S\ref{s:infinitesimal_limit}.
\end{example}

\subsection{Discrete Horizontal Lift}\label{dcpb:subsec:discrete_horizontal_lift}

The {\bfi discrete horizontal lift}\index{horizontal!lift!discrete|see{discrete horizontal lift}}\index{discrete horizontal lift} of an element $(x_0,x_1)\in S\times S$ is the subset of $Q\times Q$ that are horizontal elements, and project to $(x_0,x_1)$. Once we specify the base point $q\in Q$, the discrete horizontal lift is unique, and we introduce the map $(\cdot,\cdot)^h_q:S\times S\rightarrow Q\times Q$.

\paragraph{Discrete Horizontal Lifts from Discrete Connections.}\index{discrete horizontal lift!from discrete connection}\index{discrete connection!to discrete horizontal lift}
The discrete horizontal lift can be constructed once the discrete horizontal space is defined by a choice of discrete connection.
\begin{definition}\label{dcpb:defn:discrete_horizontal_lift}
The \textbf{discrete horizontal lift} is the unique map $(\cdot,\cdot)^h_q:S\times S\rightarrow Q\times Q$, such that
\[
(\pi,\pi)\cdot(x_0,x_1)^h_q = (x_0,x_1)\, ,\]
and
\[
(x_0,x_1)^h_q\in\operatorname{Hor}_q\, .\]
\end{definition}

\begin{lemma}
The discrete horizontal lift is $G$-equivariant, which is to say that
\[(x_0,x_1)^h_{gq} = g \cdot (x_0,x_1)^h_q \, .\]
\end{lemma}
\begin{proof}
Given $(x_0,x_1)\in S\times S$, denote $(x_0,x_1)^h_{q_0}$ by $(q_0,q_1)$. Then, by the definition of the discrete horizontal lift, we have that
\begin{align*}
(\pi,\pi)\cdot(q_0,q_1) &= (x_0,x_1)\, ,\\
\intertext{and it follows that}
(\pi,\pi)\cdot(gq_0,gq_1) &= (x_0,x_1)\, .
\end{align*}
Also, from the definition of the discrete horizontal lift,
\[
(q_0,q_1)\in\operatorname{Hor}_{q_0}\, ,
\]
and by the $G$-equivariance of the horizontal space,
\[
(gq_0,gq_1)\in g\cdot\operatorname{Hor}_{q_0} = \operatorname{Hor}_{gq_0}\, .
\]
This implies that $(gq_0,gq_1)$ satisfies the conditions for being the discrete horizontal lift of $(x_0,x_1)$ with basepoint $gq_0$. Therefore, $(x_0,x_1)^h_{gq_0}=(gq_0,gq_1)=g\cdot(q_0,q_1)=g\cdot (x_0,x_1)^h_{q_0}$, as desired.
\end{proof}

\paragraph{Discrete Connections from Discrete Horizontal Lifts.}\index{discrete connection!from discrete horizontal lift}\index{discrete horizontal lift!to discrete connection}
Conversely, given a discrete horizontal lift, we can recover a discrete connection.
\begin{definition}
Given a discrete horizontal lift, we define the \textbf{horizontal component} to be
\[ \operatorname{hor}(q_0,q_1) = (\pi(q_0,q_1))^h_{q_0}\, ,\]
and the \textbf{vertical component} is given by
\[ \operatorname{ver}(q_0,q_1) = i_{q_0}(g)\, ,\]
where $g$ is the unique group element such that
\[ (q_0,q_1) = i_{q_0}(g)\cdot \operatorname{hor}(q_0,q_1) \, .\]
\end{definition}
The last expression simply states that the discrete horizontal and vertical space are complementary with respect to the composition we defined between a vertical element and an arbitrary element of $Q\times Q$.

\paragraph{Discrete Horizontal Lifts from Discrete Connection $1$-Forms.}\index{discrete horizontal lift!from discrete connection $1$-form}\index{discrete connection $1$-form!to discrete horizontal lift}
We wish to construct a discrete horizontal lift
$(\cdot,\cdot)^h:S\times S\rightarrow (Q\times Q)/G$, given a
discrete connection $\mathcal{A}_d:Q\times Q\rightarrow G$. We state the construction
of such a discrete horizontal lift as a proposition.

\begin{proposition}
Given a discrete connection $1$-form, $\mathcal{A}_d:Q\times Q\rightarrow G$, the discrete horizontal lift is given by
\[(x_0,x_1)^h=[\pi^{-1}(x_0,x_1)\cap\mathcal{A}_d^{-1}(e)]_G.\]
Furthermore, the discrete horizontal lift satisfies the following
identity,
\[
i_{q_0}(\mathcal{A}_d(q_0,q_1))\cdot(\pi(q_0,q_1))^h_{q_0} = (q_0,q_1),
\]
which implies that the discrete connection $1$-form and the discrete horizontal lift induces a horizontal and vertical decomposition of $Q\times Q$.

The horizontal lift can be expressed in a local trivialization,
where $q_0=(x_0,g_0),$ using the local expression for the discrete
connection,
\[(x_0,x_1)^h_{q_0}=(x_0,g_0,x_1,g_0 (A(x_0,x_1))^{-1}).\]
\end{proposition}
\begin{proof}
We will show that this operation is well-defined on the
quotient space. Using the local representation of the discrete connection in the
local trivialization (see Definition~\ref{dcpb:defn:local_representation}), we have,
\begin{align*}
&\mathcal{A}_d^{-1}(e)\cap\pi^{-1}(x_0,x_1)\\
&\qquad=\{(\tilde{x}_0, g, \tilde{x}_1, g\cdot
(A(\tilde{x}_0,\tilde{x}_1))^{-1}) \mid \tilde{x}_0,\tilde{x}_1\in
S, g\in G\} \\
&\qquad\qquad\cap
\{(x_0,h_0,x_1,h_1)\mid h_0,h_1\in G\}\\
&\qquad=\{(x_0,g,x_1,g\cdot (A(x_0,x_1))^{-1})\mid h\in G\}\\
&\qquad=G\cdot(x_0,e,x_1,(A(x_0,x_1))^{-1}),
\end{align*}
which is a well-defined element of $(Q\times Q)/G$. Since this is
true in a local trivialization, and both the discrete connection
and projection operators are globally defined, this inverse coset
is globally well-defined as an element of $(Q\times Q)/G$.

In particular, the computation above allows us to obtain a local
expression for the discrete horizontal lift in terms of the local
representation of the discrete connection. That is,
\begin{align*}
(x_0,x_1)^h_{(x_0,e)} &= (x_0,e,x_1,(A(x_0,x_1))^{-1}),\\
(x_0,x_1)^h &= [(x_0,e,x_1,(A(x_0,x_1))^{-1})]_G.
\end{align*}
By the properties of the discrete horizontal lift, this extends to
$\pi^{-1}(x_0,x_1)\subset Q\times Q$,
\begin{align*}
(x_0,x_1)^h_{(x_0,g)} &= (x_0,x_1)^h_{g(x_0,e)}\\
&=g\cdot(x_0,x_1)^h_{(x_0,e)}\\
&=g(x_0,e,x_1,(A(x_0,x_1))^{-1})\\
&=(x_0,g,x_1,g(A(x_0,x_1))^{-1}).
\end{align*}

To prove the second claim, we have in the local trivialization of
$Q\times Q$, $(q_0,q_1)=(x_0,g_0,x_1,g_1)$. Then, by the result
above,
\[
(\pi(q_0,q_1))^h_{q_0}=(\pi(q_0,q_1))^h_{(x_0,g_0)}=(x_0,g_0,x_1,g_0
(A(x_0,x_1))^{-1}).
\]
Also, by the local representation of the discrete connection,
\[
\mathcal{A}_d(q_0,q_1)=g_1 A(x_0,x_1) g_0^{-1}.
\]
Therefore,
\begin{align*}
i_{q_0}(\mathcal{A}_d(q_0,q_1)\cdot (\pi(q_0,q_1))^h_{q_0}
&=(e,\mathcal{A}_d(q_0,q_1))\cdot(\pi(q_0,q_1))^h_{q_0}\\
&= (e,g_1 A(x_0,x_1) g_0^{-1})\cdot (x_0,g_0,x_1,g_0 (A(x_0,x_1))^{-1})\\
&= (x_0,g_0,x_1,(g_1 A(x_0,x_1) g_0^{-1})(g_0 (A(x_0,x_1))^{-1}))\\
&= (x_0,g_0,x_1,g_1)\\
&= (q_0,q_1),\end{align*} as claimed.
\end{proof}

\paragraph{Discrete Connection $1$-Forms from Discrete Horizontal Lifts.}\index{discrete connection $1$-form!from discrete horizontal lift}\index{discrete horizontal lift!to discrete connection $1$-form}
Given a horizontal lift $(\cdot,\cdot)^h_q:S\times
S\rightarrow Q\times Q$, we wish to construct a discrete
connection $1$-form, $\mathcal{A}_d:Q\times Q\rightarrow G$.

\begin{lemma}
Given a discrete horizontal lift, $(\cdot,\cdot)^h_q:S\times S\rightarrow Q\times Q$, the \textbf{discrete connection $1$-form}, $\mathcal{A}_d:Q\times Q\rightarrow G$, is uniquely defined by the following identity,
\[
i_{q_0}(\mathcal{A}_d(q_0,q_1))\cdot(\pi(q_0,q_1))^h_{q_0} = (q_0,q_1).
\]
\end{lemma}
\begin{proof}
To show that this construction is well-defined, we note that
$\pi_1(q_0,q_1)=\pi_1(\pi(q_0,q_1))^h_{q_0}$, by the construction
of $(\cdot,\cdot)^h_{q_0}$ from $(\cdot,\cdot)^h:S\times
S\rightarrow (Q\times Q)/G$. Furthermore, $\pi_2(q_0,q_1)$ and
$\pi_2(\pi(q_0,q_1))^h_{q_0}$ are in the same fiber of the
principal bundle $\pi:Q\rightarrow Q/G$ and are therefore related
by a unique element $g\in G$. Since this element is unique,
$\mathcal{A}_d(q_0,q_1)$ is uniquely defined by the identity.
\end{proof}

\subsection{Splitting of the Discrete Atiyah Sequence (Connection $1$-Form)}\label{dcpb:subsec:splitting_connection}\index{discrete Atiyah sequence!splitting}

Consider the discrete Atiyah sequence,
\[\xymatrix{
0\ar[r] & \tilde{G} \ar@<0.5ex>[r]^(0.32){(q,gq)}
\ar@<-0.5ex>@{<--}[r]_(0.32){(\pi_1,\mathcal{A}_d)}& (Q\times Q)/G
\ar@<0.5ex>[r]^(0.59){(\pi,\pi)} \ar@<-0.5ex>@{<--}[r]_(0.62){(\cdot,\cdot)^h}
& S \times S \ar[r] & 0}\, .\]
Given a short exact sequence
\[\xymatrix{
0\ar[r] & A_1 \ar@<0.5ex>[r]^{f} \ar@<-0.5ex>@{<--}[r]_{k} & B
\ar@<0.5ex>[r]^{g} \ar@<-0.5ex>@{<--}[r]_{h}& A_2 \ar[r] & 0 }\, ,\]
there are three equivalent conditions under which the exact sequence is split. They are as follows,
\begin{enumerate}
\item There is a homomorphism $h:A_2\rightarrow B$ with $g\circ
h=1_{A_2}$; \item There is a homomorphism $k:B\rightarrow A_1$
with $k\circ f=1_{A_1}$; \item The given sequence is isomorphic
(with identity maps on $A_1$ and $A_2$) to the direct sum short
exact sequence,
\[\xymatrix{0 \ar[r] & A_1 \ar[r]^(0.4){i_1} & A_1\oplus A_2 \ar[r]^(0.6){\pi_2} & A_2 \ar[r] & 0
},\]
and in particular, $B\cong A_1\oplus A_2$.
\end{enumerate}
We will address all three conditions in this and the next two subsections.

\paragraph{Splittings from Discrete Connection $1$-Forms.}\index{discrete Atiyah sequence!splitting!from discrete connection $1$-form}\index{discrete connection $1$-form!to splitting of Atiyah sequence}
A discrete connection $1$-form, $\mathcal{A}_d:Q\times Q\rightarrow G$, induces a splitting of the discrete Atiyah sequence, in the sense that
\[ (\pi_1,\mathcal{A}_d)\circ i = 1_{\tilde{G}}\, . \]

\begin{lemma}
Given a discrete connection $1$-form, $\mathcal{A}_d:Q\times Q\rightarrow G$, we obtain a splitting of the discrete Atiyah sequence, $\varphi:(Q\times Q)/G\rightarrow \tilde{G}$, which is given by
\[\varphi([q_0,q_1]_G)=[q_0,\mathcal{A}_d(q_0,q_1)]_G.\]
We denote this map by $(\pi_1,\mathcal{A}_d)$.
\end{lemma}
\begin{proof}
This expression is well-defined, as the following computation
shows,
\begin{align*}
\varphi([gq_0,gq_1]_G)
&=[gq_0,\mathcal{A}_d(gq_0,gq_1)]_G\\
&=[gq_0,g\mathcal{A}_d(q_0,q_1)g^{-1}]_G\\
&=\varphi([q_0,q_1]_G).
\end{align*}
Furthermore, since
\begin{align*}
(\pi_1,\mathcal{A}_d)\circ i ([q,g]_G)
&= (\pi_1,\mathcal{A}_d)([q,gq]_G)\\
&= [\pi_1(q,gq),\mathcal{A}_d(q,gq)]_G\\
&= [q,g]_G,
\end{align*}
it follows that we obtain a splitting of the discrete Atiyah sequence.
\end{proof}

\paragraph{Discrete Connection $1$-Forms from Splittings.}\index{discrete connection $1$-form!from splitting of Atiyah sequence}\index{discrete Atiyah sequence!splitting!to discrete connection $1$-form}
Given a splitting of the discrete Atiyah sequence, we can obtain a discrete connection $1$-form using the following construction.

Given $[q_0,q_1]_G\in(Q\times Q)/G$, we obtain from the splitting of the
discrete Atiyah sequence an element, $[q,g]_G\in\tilde{G}$.
Viewing $[q,g]_G$ as a subset of $Q\times G$, consider the unique
$\tilde{g}$ such that $(q_0,\tilde{g})\in [q,g]_G\subset Q\times
G$. Then, we define
\[\mathcal{A}_d(x_0,e,x_1,g_0^{-1}g_1)=\tilde{g}.\]
We extend this definition to the whole of $Q\times Q$ by
equivariance,
\[\mathcal{A}_d(x_0,g_0,x_1,g_1)=g_0\tilde{g}g_0^{-1}.\]
\begin{lemma}
Given a splitting of the discrete Atiyah sequence, the construction above yields a discrete connection $1$-form with the requisite properties.
\end{lemma}
\begin{proof}
To show that the $\mathcal{A}_d$ satisfies the properties of a
discrete connection $1$-form, we first note that equivariance follows
from the construction.

Since we have a splitting, it follows that
$\varphi([q,gq]_G)=[q,g]_G$, as $\varphi$ composed with the map
from $\tilde{G}$ to $(Q\times Q)/G$ is the identity on
$\tilde{G}$. Using a local trivialization, we have,
\begin{align*}
[q_0,g]_G &= \varphi([q_0,gq_0]_G)\\
&= \varphi([(x_0,e),(x_0,g_0^{-1} g
g_0)]_G)\\
&= [(x_0,e),\tilde{g}]_G\\
&= [(x_0,g_0),g_0 \tilde{g} g_0^{-1}]_G.
\end{align*} Then, by definition,
\[\mathcal{A}_d((x_0,e),(x_0,g_0^{-1}g g_0))=\tilde{g},\]
and furthermore, $g=g_0 \tilde{g} g_0^{-1}$. From this, we
conclude that
\begin{align*}
\mathcal{A}_d(q_0,gq_0)
&=\mathcal{A}_d((x_0,g_0),(x_0,g g_0))\\
&=g_0\mathcal{A}_d((x_0,e),(x_0,g_0^{-1} g g_0)) g_0^{-1}\\
&=g_0 \tilde{g} g_0^{-1}\\
&=g.
\end{align*}
Therefore, we have that $\mathcal{A}_d(q_0,gq_0)=g$, which
together with equivariance implies that $\mathcal{A}_d$ is a
discrete connection $1$-form.
\end{proof}

\subsection{Splitting of the Discrete Atiyah Sequence (Horizontal Lift)}\label{dcpb:subsec:splitting_lift}\index{discrete Atiyah sequence!splitting}
As was the case with the discrete connection $1$-form, the discrete horizontal lift is in one-to-one correspondence with splittings of the discrete Atiyah sequence, and they are related by taking the quotient, or extending by $G$-equivariance, as appropriate.
\paragraph{Splittings from Discrete Horizontal Lifts.}\index{discrete Atiyah sequence!splitting!from discrete horizontal lift}\index{discrete horizontal lift!to splitting of Atiyah sequence}
Given a discrete horizontal lift, we obtain a splitting by taking its quotient.
\begin{lemma}
Given a discrete horizontal lift, $(\cdot,\cdot)^h_q:S\times S\rightarrow Q\times Q$, the map $(\cdot,\cdot)^h:S\times S\rightarrow (Q\times Q)/G$, which is given by
\[ (x_0,x_1)^h = [(x_0,x_1)^h_{(x_0,e)}]_G\, ,\]
induces a splitting of the discrete Atiyah sequence.
\end{lemma}
\begin{proof}
We compute,
\begin{align*}
(\pi,\pi)\circ(x_0,x_1)^h
&= (\pi,\pi)([(x_0,x_1)^h_{(x_0,e)}]_G)\\
&= (x_0,x_1)\, ,
\end{align*}
where we used the $G$-equivariance of the discrete horizontal lift, and the property that $(\pi,\pi)\cdot(x_0,x_1)^h_q=(x_0,x_1)$ for any $q\in Q$. This implies that $(\pi,\pi)\circ(\cdot,\cdot)^h=1_{S\times S}$, as desired.
\end{proof}

\paragraph{Discrete Horizontal Lifts from Splittings.}\index{discrete horizontal lift!from splitting of Atiyah sequence}\index{discrete Atiyah sequence!splitting!to discrete horizontal lift}
Given a splitting, $(\cdot,\cdot)^h:S\times S\rightarrow (Q\times Q)/G$, we obtain a discrete horizontal lift, $(\cdot,\cdot)^h_q:S\times S\rightarrow Q\times Q$, using the following construction.

We denote by $(x_0,x_1)^h_{q_0}$ the unique element in $(x_0,x_1)^h$, thought of as a subset of $Q\times Q$, such that the first component is $q_0$. This is the discrete horizontal lift of the point
$(x_0,x_1)\in S\times S$ where the base point is specified.

\begin{lemma}
Given a splitting of the discrete Atiyah sequence, the construction above yields a discrete horizontal lift with the requisite properties.
\end{lemma}
\begin{proof}
Since the quotient space $(Q\times Q)/G$ is obtained by the diagonal action of $G$ on $Q\times Q$, it follows that if $(x_0,x_1)^h_{q_0}\in (x_0,x_1)^h\subset Q\times Q$, then $g\cdot(x_0,x_1)^h_{q_0}\in (x_0,x_1)^h\subset Q\times Q$. Since the first component of $G\cdot(x_0,x_1)^h_{q_0}$ is $gq_0$, and  $g\cdot(x_0,x_1)^h_{q_0}\in (x_0,x_1)^h\subset Q\times Q$, we have that
\[ (x_0,x_1)^h_{gq_0}= g\cdot(x_0,x_1)^h_{q_0}\, ,\]
which is to say that the discrete horizontal lift constructed above is $G$-equivariant.

Since $(\cdot,\cdot)^h$ is a splitting of the discrete Atiyah sequence, we have that $(\pi,\pi)\circ(\cdot,\cdot)^h=1_{S\times S}$, and this implies that any element in $(x_0,x_1)^h$, viewed as a subset of $Q\times Q$, projects to $(x_0,x_1)$. Therefore, the discrete horizontal lift we constructed above has the requisite properties.
\end{proof}

\subsection{Isomorphism between $(Q\times Q)/G$ and $(S\times S)\oplus\tilde{G}$}\label{dcpb:subsec:isomorphism}\index{discrete connection!isomorphism}
The notion of a discrete connection is motivated by the desire to
construct a global diffeomorphism between $(Q\times
Q)/G\rightarrow S$ and $(S\times S)\oplus\tilde{G}\rightarrow S.$
This is the discrete analogue of the identification between
$TQ/G\rightarrow Q/G$ and
$T(Q/G)\oplus\tilde{\mathfrak{g}}\rightarrow Q/G$ which is the
context for Lagrangian Reduction in \cite{CeMaRa2001}. Since a choice of discrete connection corresponds to a choice of splitting of the discrete Atiyah sequence, we have the following commutative diagram, where each row is a short exact sequence.
\[\xymatrix{
0\ar[r] & \tilde{G} \ar@{=}[d]_{1_{\tilde{G}}}
\ar@<0.5ex>[r]^(0.32){(q,gq)}
\ar@<-0.5ex>@{<--}[r]_(0.32){(\pi_1,\mathcal{A}_d)}& (Q\times Q)/G
\ar@<0.5ex>[r]^(0.59){(\pi,\pi)}
\ar@<-0.5ex>@{<--}[r]_(0.62){(\cdot,\cdot)^h}
\ar[d]^{\alpha_{\mathcal{A}_d}}& S \times S \ar@{=}[d]^{1_{S
\times S}}
\ar[r] & 0\\
0\ar[r] & \tilde{G} \ar@<0.5ex>[r]^(0.3){i_1}
\ar@<-0.5ex>@{<--}[r]_(0.3){\pi_1} & \tilde{G}\oplus(S\times S)
\ar@<0.5ex>[r]^(0.6){\pi_2} \ar@<-0.5ex>@{<--}[r]_(0.6){i_2}& S
\times S \ar[r] & 0 }\]
Here, we see how the identification between $(Q\times Q)/G$ and $(S\times S)\oplus\tilde{G}$ are
naturally related to the discrete connection and the discrete horizontal lift.

Recall that the discrete adjoint bundle $\tilde{G}$ is the associated bundle one obtains when $M=G$, and $\rho_g$ acts by conjugation. Furthermore, the action of $G$ on $Q\times Q$ is by the diagonal action, and the action of $G\times G$ on $Q\times Q$ is component-wise.

\begin{proposition}\label{dcpb:prop:isomorphism}
The map $\alpha_{\mathcal{A}_d}:(Q\times Q)/G\rightarrow (S\times
S)\oplus \tilde{G}$ defined by
\[
\alpha_{\mathcal{A}_d}([q_0,q_1]_G)=(\pi q_0, \pi
q_1)\oplus[q_0,\mathcal{A}_d(q_0,q_1)]_G,
\] is a well-defined bundle isomorphism.
The inverse of $\alpha_{A_d}$ is given by
\[
\alpha_{\mathcal{A}_d}^{-1}((x_0,x_1)\oplus[q,g]_G)=[(e,g)\cdot(x_0,x_1)^h_q
]_G,
\] for any
$q\in Q$ such that $\pi q=x_0$.
\end{proposition}
\begin{proof}
To show that $\alpha_{\mathcal{A}_d}$ is well-defined, note that
for any $g\in G$, we have that
\[
(\pi gq_0, \pi gq_1)=(\pi q_0, \pi q_1),
\]
and also,
\begin{align*}
[gq_0,\mathcal{A}_d(gq_0,gq_1)]_G &= [gq_0,g\mathcal{A}_d(q_0,q_1)g^{-1}]_G\\
&=[q_0,\mathcal{A}_d(q_0,q_1)]_G.
\end{align*}
Then, we see that
\[ \alpha_{\mathcal{A}_d}([gq_0,gq_1]_G) =\alpha_{\mathcal{A}_d}([q_0,q_1]_G). \]
To show that $\alpha_{\mathcal{A}_d}^{-1}$ is well-defined, note
that for any $k\in G$,
\[
(x_0,x_1)^h_{kq}=k\cdot(x_0,x_1)^h_{q},
\]
and that
\begin{align*}
\alpha_{\mathcal{A}_d}^{-1}((x_0,x_1)\oplus[kq,kgk^{-1}]_G)
& = [(e,kgk^{-1})\cdot(x_0,x_1)^h_{kq}]_G\\
& = [(e,kgk^{-1})\cdot k\cdot(x_0,x_1)^h_q]_G\\
& = [(ek,kgk^{-1}k)\cdot (x_0,x_1)^h_q]_G\\
& = [(ke, kg)\cdot(x_0,x_1)^h_q]_G\\
& = [k\cdot(e,g)\cdot(x_0,x_1)^h_q]_G\\
& = [(e,g)\cdot(x_0,x_1)^h_q]_G\\
& = \alpha_{\mathcal{A}_d}^{-1}((x_0,x_1)\oplus[q,g]_G).\qedhere
\end{align*}
\end{proof}

\begin{example}
It is illustrative to consider the notion of a discrete connection,
and the isomorphism, in the degenerate case when $Q=G$, which is
the context of discrete Euler--Poincar\'{e} reduction. Here, the
isomorphism is between $(G\times G)/G$ and $\tilde G$, and the connection
$\mathcal{A}_d:G\times G\rightarrow G$ is given by
\[
\mathcal{A}_d(g_0,g_1)=g_1\cdot g_0^{-1}.
\]
Then, we have that
\begin{align*}
\alpha_{\mathcal{A}_d}([g_0,g_1]_G)
&= (\pi g_0, \pi g_1)\oplus[g_0,\mathcal{A}_d(q_0,q_1)]_G\\
&= (e,e)\oplus [g_0, g_1 g_0^{-1}]_G\, .
\end{align*}
Taking the inverse, we have,
\begin{align*}
\alpha_{\mathcal{A}_d}^{-1}([g_0,g_1 g_0^{-1}]_G)
&= [(e,g_1 g_0^{-1}]\cdot(e,e)^h_{g_0}]_G\\
&= [(e,g_1 g_0^{-1}]\cdot(g_0,g_0)]_G\\
&= [e g_0, g_1 g_0^{-1} g_0 ]_G\\
&= [ g_0, g_1]_G\, ,
\end{align*}
as expected.
\end{example}

\subsection{Discrete Horizontal and Vertical Subspaces Revisited}
Having now fully introduces all the equivalent representations of a discrete connection, we can revisit the notion of discrete horizontal and vertical subspaces in light of the new structures we have introduced.

Consider the following split exact sequence,
\[\xymatrix{
0\ar[r] & A_1 \ar@<0.5ex>[r]^{f} \ar@<-0.5ex>@{<--}[r]_{k} & B
\ar@<0.5ex>[r]^{g} \ar@<-0.5ex>@{<--}[r]_{h}& A_2 \ar[r] & 0 }\, .\]
We can decompose any element in $B$ into a $A_1$ and $A_2$ term by
considering the following isomorphism,
\[B\cong f\circ k (B) \oplus h\circ g(B).\]
Similarly, in the discrete Atiyah sequence, we can decompose an
element of $(Q\times Q)/G$ into a horizontal and vertical piece by
performing the analogous construction on the split exact sequence
\[\xymatrix{
0\ar[r] & \tilde{G} \ar@<0.5ex>[r]^(0.32){(q,gq)}
\ar@<-0.5ex>@{<--}[r]_(0.32){(\pi_1,\mathcal{A}_d)}& (Q\times Q)/G
\ar@<0.5ex>[r]^(0.59){(\pi,\pi)}
\ar@<-0.5ex>@{<--}[r]_(0.62){(\cdot,\cdot)^h} & S \times S \ar[r]
& 0}\, .\]
This allows us to define horizontal and vertical spaces associated with the pair groupoid $Q\times Q$, in terms of all the structures we have introduced.
\begin{definition}
The \textbf{horizontal space} is given by
\begin{align*}
\operatorname{Hor}_q
&=\{(q,q')\in Q\times Q \mid \mathcal{A}_d(q,q')=e\}\\
&=\{(\pi q,x_1)^h_q\in Q\times Q \mid x_1\in S\}.
\end{align*}
This is the discrete analogue of the statement $\operatorname{Hor}_q=\{v_q\in TQ\mid
\mathcal{A}(v_q)=0\}=\{(v_{\pi q})^h_q\in TQ\mid v_{\pi q}\in TS\}$.
\end{definition}
\begin{definition}
The \textbf{vertical space} is given by
\begin{align*} \operatorname{Ver}_q
&=\{(q,q')\in Q\times Q \mid
(\pi,\pi)(q,q')=e_{S\times S}\}\\
&=\{i_q(g)\mid g\in G\}.
\end{align*}
This is the discrete analogue of the statement
$\operatorname{Ver}_q=\{v_q\in TQ\mid \pi_*(v_q)=0\}=\{\xi_Q\mid
\xi\in\mathfrak{g}\}$.
\end{definition}
In particular, we can decompose an element of $Q\times Q$ into a
horizontal and vertical component.
\begin{definition}
The \textbf{horizontal component}\index{discrete connection!horizontal component}\index{horizontal!component!discrete|see{discrete connection, horizontal component}} of $(q_0,q_1)\in Q\times Q$ is
given by
\[\operatorname{hor}(q_0,q_1)=((\cdot,\cdot)^h\circ(\pi,\pi))(q_0,q_1)=(\pi q_0,\pi q_1)^h_{q_0}.\]
\end{definition}
\begin{definition}
The \textbf{vertical component}\index{discrete connection!vertical component}\index{vertical!component!discrete|see{discrete connection, vertical component}} of $(q_0,q_1)\in Q\times Q$ is
given by
\[\operatorname{ver}(q_0,q_1)=(i\circ(\pi_1,\mathcal{A}_d))(q_0,q_1)=(q_0,\mathcal{A}_d(q_0,q_1)q_0)=i_{q_0}(\mathcal{A}_d(q_0,q_1)).\]
\end{definition}

\begin{lemma}
The horizontal component can be expressed as
\[\operatorname{hor}(q_0,q_1)=
i_{q_0}((\mathcal{A}_d(q_0,q_1))^{-1})\cdot (q_0,q_1).\]
\end{lemma}
\begin{proof}
\begin{align*}
i_{q_0}(\mathcal{A}_d(q_0,q_1)^{-1})\cdot(q_0,q_1) &=
(q_0,(\mathcal{A}_d(q_0,q_1))^{-1}q_0)\cdot(q_0,q_1)\\
&= (e,(\mathcal{A}_d(q_0,q_1))^{-1})(q_0,q_1)\\
&= (q_0,(\mathcal{A}_d(q_0,q_1))^{-1}q_1).
\end{align*}
Clearly,
$(\pi,\pi)(q_0,(\mathcal{A}_d(q_0,q_1))^{-1}q_1)=(\pi,\pi)(q_0,q_1)$.
Furthermore,
\[\mathcal{A}_d(q_0,(\mathcal{A}_d(q_0,q_1))^{-1}q_1)=(\mathcal{A}_d(q_0,q_1))^{-1}\mathcal{A}_d(q_0,q_1)=e.\]
Therefore, by definition,
$(q_0,(\mathcal{A}_d(q_0,q_1))^{-1}q_1)=\operatorname{hor}(q_0,q_1)$.
\end{proof}

\begin{lemma}
The horizontal and vertical operators satisfy the following
identity,
\[\operatorname{ver}(q_0,q_1)\cdot\operatorname{hor}(q_0,q_1)=(q_0,q_1).\]
\end{lemma}
\begin{proof}
\begin{align*}
\operatorname{ver}(q_0,q_1)\cdot\operatorname{hor}(q_0,q_1) &=
 i_{q_0}(\mathcal{A}_d(q_0,q_1))\cdot(i_{q_0}((\mathcal{A}_d(q_0,q_1))^{-1})\cdot(q_0,q_1))\\
&=(e,\mathcal{A}_d(q_0,q_1))(e,(\mathcal{A}_d(q_0,q_1))^{-1})(q_0,q_1)\\
&=(e,\mathcal{A}_d(q_0,q_1))(q_0,(\mathcal{A}_d(q_0,q_1))^{-1}q_1)\\
&=(q_0,\mathcal{A}_d(q_0,q_1)(\mathcal{A}_d(q_0,q_1))^{-1}q_1)\\
&=(q_0,q_1),
\end{align*}
as desired.
\end{proof}


\section{Geometric Structures Derived from the Discrete Connection}\index{discrete connection!derived geometric structures}
In this section, we will introduce some of the additional geometric structures that can be derived from a choice of discrete connection. These structures include an extension of the pair groupoid composition to take into account the principal bundle structure, continuous connections that are a limit of a discrete connection, and higher-order connection-like structures.

\subsection{Extending the Pair Groupoid Composition}\index{groupoid!composition!extended|see{discrete connection, extended groupoid composition}}\index{discrete connection!extended groupoid composition}
Recall that the composition of a vertical element $(q_0,gq_0)$ with an element $(q_0,q_1)$ is given by
\[(q_0,gq_0)\cdot(q_0,q_1)=(q_0,gq_1).\]
The choice of a discrete connection allows us to further extend the composition, in a manner that is relevant in describing the curvature of a discrete connection. The decomposition of an element of $Q\times Q$ into a horizontal and vertical piece naturally suggests a generalization of the composition operation on $Q\times Q$ (viewed as a pair groupoid), by using the discrete connection, and the principal bundle structure of $Q$.

We wish to define a composition on $Q\times Q$ such that the composition of $(q_0,q_1)\cdot(\tilde{q}_0,\tilde{q}_1)$ is defined whenever $\pi q_1=\pi\tilde{q}_0$. Furthermore, we require that the extended composition be consistent with the vertical composition
we introduced previously, as well as the pair groupoid composition, whenever their domains of definition coincide.

The extended composition is obtained by left translating $(\tilde{q}_0,\tilde{q}_1)$ by a group element $h$, such that $q_1=h\tilde{q}_0$, and then using the pair groupoid composition on $(q_0,q_1)$ and the left translated term $h(\tilde{q}_0, \tilde{g}_1)$. This yields the following intrinsic definition of the extended composition.

\begin{definition}
The \textbf{extended pair groupoid composition} of $(q_0,q_1),(\tilde{q}_0,\tilde{q}_1)\in Q\times Q$ is defined whenever $\pi q_1=\pi \tilde{q}_0$, and it is given by
\[ (q_0,q_1)\cdot (\tilde{q}_0,\tilde{q}_1) = (q_0, \mathcal{A}_d(\tilde{q}_0,q_1) \tilde{q}_1).\]
\end{definition}

As the following lemmas show, this extended composition is consistent with the vertical composition and the pair groupoid composition.

\begin{lemma}
The extended pair groupoid composition is consistent with the composition of a vertical element with an arbitrary element.
\end{lemma}
\begin{proof}
Consider the composition of a vertical element with an arbitrary element,
\begin{align*}
(q_0,gq_0)\cdot(q_0,q_1) &= (q_0,gq_1).\\
\intertext{This is consistent with the result using the extended composition,}
(q_0,gq_0)\cdot(q_0,q_1) &= (q_0, \mathcal{A}_d(q_0,gq_0)q_1)\\*
&= (q_0, g q_1)\, ,
\end{align*}
where we used that the discrete connection yields a splitting of the Atiyah sequence.
\end{proof}

\begin{lemma}
The extended pair groupoid composition is consistent with the pair groupoid composition.
\end{lemma}
\begin{proof}
The pair groupoid composition is given by
\begin{align*}
(q_0,q_1)\cdot(q_1,q_2) &= (q_0, q_2)\, .\\
\intertext{This is consistent with the extended composition,}
(q_0,q_1)\cdot(q_1,q_2) &= (q_0, \mathcal{A}_d(q_1,q_1) q_2)\\
&= (q_0, e q_2)\\
&= (q_0, q_2)\, .\qedhere
\end{align*}
\end{proof}

The extended composition is $G$-equivariant, and is well-defined on the quotient space, as the following lemma shows.

\begin{lemma}
The composition $\cdot:(Q\times Q)\times(Q\times Q)\rightarrow
(Q\times Q)$ is $G$-equivariant, that is,
\[(gq_0,gq_1)\cdot(g\tilde{q}_0,g\tilde{q}_1)=g\cdot((q_0,q_1)\cdot (\tilde{q}_0,\tilde{q}_1)).\]
Furthermore, the composition induces a well-defined quotient
composition $\cdot:((Q\times Q)\times(Q\times Q))/G\rightarrow
(Q\times Q)/G$.
\end{lemma}
\begin{proof}
Given $g\in G$, we consider,
\begin{align*}
(gq_0,gq_1)\cdot(g\tilde{q}_0,g\tilde{q}_1)
&= (g q_0, \mathcal{A}_d(g\tilde{q}_0, gq_1)g\tilde{q}_1)\\
&= (g q_0, g\mathcal{A}_d(\tilde{q}_0,q_1)g^{-1}g\tilde{q}_1)\\
&= (g q_0, g \mathcal{A}_d(\tilde{q}_0,q_1)\tilde{q}_1)\\
&= g \cdot (q_0, \mathcal{A}_d(\tilde{q}_0,q_1)\tilde{q}_1)\\
&= g \cdot ((q_0,q_1)\cdot (\tilde{q}_0, \tilde{q}_1))\, ,
\end{align*}
where we used the equivariance of the discrete connection. It
follows that the composition is equivariant. Furthermore,
\[
[(g\tilde{q}_0,g\tilde{q}_1)\cdot(gq_0,gq_1)]_G=[(\tilde{q}_0,\tilde{q}_1)\cdot(q_0,q_1)]_G,
\]
which means that $\cdot:((Q\times Q)\times(Q\times
Q))/G\rightarrow (Q\times Q)/G$ is well-defined.
\end{proof}

\begin{corollary}
The composition of $n$-terms is $G$-equivariant. That is to say,
\begin{multline*}
(gq^1_0,gq^1_1)\cdot(gq^2_0,gq^2_1)\cdot\ldots\cdot(gq^{n-1}_0,gq^{n-1}_1)\cdot(gq^n_0,gq^n_1)\\
=g\cdot((q^1_0,q^1_1)\cdot(q^2_0,q^2_1)\cdot\ldots\cdot(q^{n-1}_0,q^{n-1}_1)\cdot(q^n_0,q^n_1)).
\end{multline*}
\end{corollary}
\begin{proof}
The result follows by induction on the previous lemma.
\end{proof}

We find that the extended composition we have constructed on the pair groupoid is associative. However, as we shall see in \S\ref{dcpb:subsec:levi_civita}, composing pair groupoid elements about a loop in the shape space will not in general yield the identity element $e_{Q\times Q}$, and the defect represents the holonomy about the loop, which is related to curvature. This may yield the discrete analogue of the expression giving the geometric phase in terms of a loop integral (in shape space) of the curvature of the connection.

\begin{lemma}
The composition $\cdot:(Q\times Q)\times(Q\times Q)\rightarrow (Q\times Q)$ is
associative. That is,
\[((q_0^0,q_1^0)\cdot(q_0^1,q_1^1))\cdot(q_0^2,q_1^2)=(q_0^0,q_1^0)\cdot((q_0^1,q_1^1)\cdot(q_0^2,q_1^2))\, .\]
\end{lemma}
\begin{proof}
Evaluating the left-hand side, we obtain
\begin{align*}
((q_0^0,q_1^0)\cdot(q_0^1,q_1^1))\cdot(q_0^2,q_1^2)
&= (q_0^0,\mathcal{A}_d(q_0^1,q_1^0) q_1^1)\cdot(q_0^2, q_1^2)\\
&= (q_0^0, \mathcal{A}_d(q_0^2, \mathcal{A}_d(q_0^1, q_1^0) q_1^1) q_1^2)\\
&= (q_0^0, \mathcal{A}_d(q_0^1,q_1^0)\mathcal{A}_d(q_0^2, q_1^1)q_1^2)\, ,\\
\intertext{and the right-hand side is given by}
(q_0^0,q_1^0)\cdot((q_0^1,q_1^1)\cdot(q_0^2,q_1^2))
&= (q_0^0,q_1^0)\cdot (q_0^1, \mathcal{A}_d(q_0^2,q_1^1)q_1^2)\\
&= (q_0^0, \mathcal{A}_d(q_0^1,q_1^0)\mathcal{A}_d(q_0^2,q_1^1)q_1^2)\, .\\
\intertext{Therefore,}
((q_0^0,q_1^0)\cdot(q_0^1,q_1^1))\cdot(q_0^2,q_1^2)&=(q_0^0,q_1^0)\cdot((q_0^1,q_1^1)\cdot(q_0^2,q_1^2))\, ,
\end{align*}
and the extended groupoid composition is associative.
\end{proof}


\subsection{Continuous Connections from Discrete
Connections}\label{s:infinitesimal_limit}\index{discrete connection!to continuous connection}\index{connection!from discrete connection} Given a discrete
$G$-valued connection $1$-form, $\mathcal{A}_d:Q\times Q\rightarrow
G$, we associate to it a continuous $\mathfrak{g}$-valued
connection $1$-form, $\mathcal{A}:TQ\rightarrow\mathfrak{g}$, by the
following construction,
\[\mathcal{A}([q(\cdot)])=[\mathcal{A}_d(q(0),q(\cdot))],\] where
$[\cdot]$ denotes the equivalence class of curves associated with
a tangent vector.

This uses the intrinsic representation of the tangent bundle, which is obtained by identifying a tangent vector at a point on the manifold with the equivalence class of curves on the manifold going through the point, such that the tangent to the curve at the point is given by the tangent vector, which was illustrated earlier in Figure~\ref{dcpb:fig:curvesvectors} on page~\pageref{dcpb:fig:curvesvectors}.

More explicitly, given $v_q\in TQ$, we consider an associated
curve $q:[0,1]\rightarrow Q$, and construct the curve
$g:[0,1]\rightarrow G$, given by
\[g(t)=\mathcal{A}_d(q(0),q(t)).\]
Then,
\[\mathcal{A}(v_q)=\left.\frac{d}{dt}\right|_{t=0} g(t).\]

When computing the equations in discrete reduction theory, it is often necessary to consider horizontal and vertical variations, which we introduce below.
\begin{definition}
We introduce the \textbf{vertical variation} of a point $(q_0,q_1)\in
Q\times Q$.  Given a curve $q_1^{\epsilon}:[0,1]\rightarrow Q$, such that
$q_1^{\epsilon}(0)=q_1$, the vertical variation is given by
\[
\operatorname{ver}
\delta{q}=\left.\frac{d}{d\epsilon}\right|_{\epsilon=0}\operatorname{ver}(q_0,q_1^\epsilon)
= \left.\frac{d}{d\epsilon}\right|_{\epsilon=0} i_{q_0}(\mathcal{A}_d(q_0,q_1^\epsilon))\, .
\]
\end{definition}
\begin{definition}
We introduce the \textbf{horizontal variation} of a point $(q_0,q_1)\in
Q\times Q$.  Given a curve $q_1^{\epsilon}:[0,1]\rightarrow Q$, such that
$q_1^{\epsilon}(0)=q_1$, the horizontal variation is given by
\[
\operatorname{hor}
\delta{q}=\left.\frac{d}{d\epsilon}\right|_{\epsilon=0}\operatorname{hor}(q_0,q_1^\epsilon)
= \left.\frac{d}{d\epsilon}\right|_{\epsilon=0} (\pi(q_0,q_1^\epsilon))^h_{q_0}\, .
\]
\end{definition}

\subsection{Connection-Like Structures on Higher-Order Tangent Bundles}\index{discrete connection!higher-order tangent bundle}

Given a continuous connection, we can construct connection-like
structures on higher-order tangent bundles. This construction is
described in detail in Lemma 3.2.1 of \cite{CeMaRa2001}. In particular, given a connection $1$-form, $\mathcal{A}:TQ\rightarrow \mathfrak{g}$, we obtain a well-defined
map, $\mathcal{A}^k:T^{(k)}Q\rightarrow k\mathfrak{g}$.

As we will see later, these connection-like structures on higher-order tangent bundles will provide an intrinsic method of characterizing the order of approximation of a continuous
connection by a discrete connection.

We will describe the discrete analogue of this construction. To
begin, the discrete analogue of the $k$-th order tangent bundle,
$T^{(k)}Q$, is $k+1$ copies of $Q$, namely $Q^{k+1}$. Intermediate
spaces between $T^{(k)}Q$ and $Q^{k+1}$ arise in the general
theory of multi-spaces, which is introduced in~\cite{Olver2001}.

The discrete analogue of tangent lifts, and their higher-order
analogues, are obtained by componentwise application of the map,
since a tangent lift of a map is computed by applying the map to a
representative curve, and taking its equivalence class. Therefore,
given a map $f:M\rightarrow N$, we have the naturally induced map,
\[T^{(k)}f:M^{k+1}\rightarrow N^{k+1}\quad \text{given by} \quad
T^{(k)}f(m_0,\ldots,m_k)=(f(m_0),\ldots,f(m_k)).\]
And in particular, the group action is lifted to the diagonal group action on the product space.

The discrete connection can be extended to $Q^{k+1}$ in the
natural way, $\mathcal{A}_d^k:Q^{k+1}\rightarrow
\oplus_{l=0}^{k-1} G\equiv kG$,
\[\mathcal{A}_d^k(q_0,\ldots,q_k)=\oplus_{l=0}^{k-1}
\mathcal{A}_d(q_l,q_{l+1}).\]

Similarly, we can define the map from $Q^{k+1}$ to the Whitney sum
of $k$ copies of the conjugate bundle $\tilde{G}$ by
\[Q^{k+1}\rightarrow k\tilde{G}\quad\text{by}\quad
(q_0,\ldots,q_k)\mapsto\oplus_{l=0}^{k-1}[q_0,\mathcal{A}_d(q_l,q_{l+1})]_G.\]

In a natural way, we have the following proposition.
\begin{proposition}\label{dcpb:prop:higher_order_isomorphism}
The map
\[\alpha_{\mathcal{A}_d^k}:Q^{k+1}\rightarrow
(Q/G)^{k+1}\times_{Q/G}k\tilde{G}\] defined by
\[\alpha_{\mathcal{A}_d^k}(q_0,\ldots,q_k)=(\pi q_0,\ldots, \pi
q_k)\times_{Q/G}\oplus_{l=0}^{k-1}[q_0,\mathcal{A}_d(q_l,q_{l+1})]_G,\]
is a well-defined bundle isomorphism. The inverse of
$\alpha_{\mathcal{A}_d^k}$ is given by
\begin{align*}&\alpha_{\mathcal{A}_d^k}^{-1}\left((x_0,\ldots,x_k)\times_{Q/G}\oplus_{l=0}^{k-1}\left[q_l,g_l\right]_G\right)\\
&\qquad=[(e,g_0,g_1 g_0,\ldots, g_{k-1}\ldots
g_0))\cdot(x_0,\ldots,x_k)^h_{q_0}]_G,
\end{align*}
where $(x_0,\ldots,x_k)^h_{q_0}=(\bar{q}_0,\ldots,\bar{q}_k)$ is
defined by the conditions:
\begin{align*}
\bar{q}_0 &= q_0,\\
\pi\bar{q}_l &= x_l,\\
\mathcal{A}_d(\bar{q}_l,\bar{q}_{l+1})&=e.
\end{align*}
\end{proposition}
\begin{remark}
In a local trivialization, where $q_0=(h_0,x_0)$, we have,
\[\bar{q}_{l+1}=\left((A(x_l,x_{l+1}) \cdots A(x_0,x_1)
)^{-1}h_0,x_{l+1}\right).\]
\end{remark}


\section{Computational Aspects}\index{discrete connection!computation}
While we saw in the previous section that a discrete connection induces a continuous connection in the limit, we are often concerned with constructing discrete connections that approximate a continuous connection to a given order of approximation. This section will address the question of what it means for a discrete connection to approximate a continuous connection to a given order, as well as introduce methods for constructing such discrete connections.

\subsection{Exact Discrete
Connection}\index{discrete connection!exact} It is interesting from the point of view of
computation to construct an {\bfi exact discrete connection}
associated with a prescribed continuous connection, so that we can
make sense of the statement that a given discrete connection is a
$k$-th order approximation of a continuous connection.

\paragraph{Additional Structure.}
To construct the exact discrete connection, we require that the configuration manifold $Q$ be
endowed with a bi-invariant Riemannian metric, with the property
that the associated exponential map,
\[\exp:TQ\rightarrow Q,\]
is consistent with the group action, in the sense that
\[\exp(\xi_Q(q))=\exp(\xi)\cdot q.\]
We extend the exponential to $Q\times Q$ as follows,
\begin{align*}
\overline{\exp}:TQ & \rightarrow Q\times Q,\\
v_q & \mapsto(q,\exp(v_q)),
\end{align*}
and denote the inverse by $\overline{\log}:Q\times Q\rightarrow
TQ$, which is defined in a neighborhood of the diagonal of $Q\times Q$.

\paragraph{Construction.}
Having introduced the appropriate structure on the configuration manifold, we define the exact discrete connection as follows.

\begin{definition}
The \textbf{Exact Discrete Connection} $\mathcal{A}_d^E$
associated with a prescribed continuous connection
$\mathcal{A}:TQ\rightarrow\mathfrak{g}$ and a Riemannian metric is given by
\[\mathcal{A}_d^E(q_0,q_1)=\exp(\mathcal{A}(\overline{\log}(q_0,q_1))).\]
This construction is more clearly illustrated in the following
diagram,
\[\xymatrix{
Q\times Q \ar[r]_(0.55){\overline{\log}} \ar`u[r]
`[rrr]^(0.59){\mathcal{A}_d^E} [rrr] & TQ \ar[r]_{\mathcal{A}} &
\mathfrak{g} \ar[r]_{\exp} & G }\]
\end{definition}

\paragraph{Properties.}
The exact discrete connection satisfies the properties of the discrete connection $1$-form, in that it is $G$-equivariant, and it induces a splitting of the discrete Atiyah sequence. The equivariance of the exact discrete connection arises from the fact that each of the composed maps is equivariant, and the splitting condition,
\[\mathcal{A}_d^E(i_q(g))=g,\]
is a consequence of the compatibility condition,
\[\exp(\xi_Q(q))=\exp(\xi)\cdot q.\]
Since the logarithm map is only defined on a neighborhood of the diagonal of $Q\times Q$, it follows that the exact discrete connection will have the same restriction on its domain of definition.

\begin{example}[Discrete Mechanical Connection]\label{ex:discrete_mechanical_connection}\index{discrete connection!mechanical} The continuous mechanical
connection is defined by the following diagram,
\[\xymatrix{
T^{*}Q \ar[r]^{J} & \mathfrak{g}^{*}\\
TQ \ar[r]_{\mathcal{A}} \ar[u]^{\mathbb{F}L} & \mathfrak{g}
\ar[u]_{\mathbb{I}}}\] Correspondingly, the discrete mechanical
connection is defined by the following diagram,
\[\xymatrix{
Q\times Q \ar[r]_(0.53){\mathbb{F}L_d} \ar@{-<}`u[r]
`[rr]^(0.53){J_d} [rr] & T^*Q \ar[r]_(0.52){J} &
\mathfrak{g}^* \\
Q\times Q \ar@{->}`d[r] `[rrr]_(0.53){\mathcal{A}_d} [rrr]
\ar@{=}[u] \ar[rr] & & \mathfrak{g} \ar[r]^{\exp}
\ar[u]^{\mathbb{I}} & G }\]
This is consistent with our notion of an exact discrete mechanical
connection as the following diagram illustrates,
\[\xymatrix{
Q\times Q \ar[r]_(0.53){\mathbb{F}L_d} \ar@{-<}`u[r]
`[rr]^(0.53){J_d} [rr] & T^*Q \ar[r]_(0.52){J} &
\mathfrak{g}^* \\
Q\times Q \ar@{->}`d[r] `[rrr]_(0.53){\mathcal{A}_d} [rrr]
\ar@{=}[u] \ar[r]^{\overline{\log}} & TQ \ar[u]_{\mathbb{F}L}
\ar[r]^{\mathcal{A}} & \mathfrak{g} \ar[r]^{\exp}
\ar[u]^{\mathbb{I}} & G \save "1,2"."2,3"*[F.]\frm{} \restore }\]
where the portion in the dotted box recovers the continuous
mechanical connection. In checking $G$-equivariance, we use the equivariance of
$\exp:\mathfrak{g}\rightarrow G,$ $J_d:Q\times Q\rightarrow
\mathfrak{g}^*,$ and the equivariance of $\mathbb{I}$ in the sense
of a map $\mathbb{I}:Q\rightarrow L(\mathfrak{g},\mathfrak{g}^*),$
namely, $\mathbb{I}(gq)\cdot
\operatorname{Ad}_g\xi=\operatorname{Ad}^*_{g^{-1}}\mathbb{I}(q)\cdot\xi.$
\end{example}

\subsection{Order of Approximation of a Connection}\index{discrete connection!order of approximation} We have the
necessary constructions to consider the order to which a discrete
connection approximates a continuous connection. There are two
equivalent ways of defining the order of approximation of a
continuous connection by a discrete connection, the first is more
analytical, and is given by the order of convergence in an
appropriate norm on the group.

\begin{definition}[Order of Connection, Analytic]
A discrete connection $\mathcal{A}_d$ is a $k$-th order discrete
connection if, $k$ is the maximum integer for which
\[\exists 0<c<\infty,\]
\[\exists h_0>0,\]
such that
\[\sup_{\scriptsize\begin{matrix}
  v_q\in TQ, \\
  |v_q|=1
\end{matrix}} \| \mathcal{A}_d^E(q,\exp(h v_q))(\mathcal{A}_d(q,\exp(h
v_q)))^{-1} \|\leq c h^{k+1},\quad \forall h<h_0.\]
\end{definition}

The second definition is more intrinsic, and is related to
considering the infinitesimal limit of a discrete connection to
connection-like structures on higher-order tangent bundles,
without the need for the introduction of the exact discrete
connection.

Recall from \S\ref{s:infinitesimal_limit} that we can construct a
continuous connection from a discrete connection by the following
construction,
\[\mathcal{A}([q(\cdot)])=[\mathcal{A}_d(q(0),q(\cdot))].\]
Given $\mathcal{A}_d^k:Q^{k+1}\rightarrow kG$, we can obtain
the continuous limit $\mathcal{A}^k:T^{(k)}Q\rightarrow
k\mathfrak{g}$ in a similar fashion.

\begin{definition}[Order of Connection, Intrinsic] A discrete
connection $\mathcal{A}_d$ is a $k$-th order approximation to
$\mathcal{A}$ if, $k$ is the maximum integer for which the diagram
holds,
\[\xymatrix@!0@R=1.5cm@C=5cm{
\mathcal{A}_d:Q\times Q\rightarrow G \ar@{=>}[d] &
\mathcal{A}:TQ\rightarrow \mathfrak{g} \ar@{=>}[d]\\
\mathcal{A}_d^k:Q^{k+1}\rightarrow kG \ar@{-->}[r]&
\mathcal{A}^k:T^{(k)}Q\rightarrow k\mathfrak{g} }\]
Here, the double arrows represent the higher-order structures induced by the connections, and the dotted arrow represents convergence in the limit.
\end{definition}

\subsection{Discrete Connections from Exponentiated Continuous Connections}\label{dcpb:subsec:exponentiated_discrete_connections}\index{discrete connection!from continuous connection}\index{connection!to discrete connection}
To apply the exponential and logarithm approach to construct a discrete connection from a prescribed connection, in the sense of the diagram,
\[\xymatrix{
Q\times Q \ar[r]_(0.55){\overline{\log}} \ar`u[r]
`[rrr]^(0.59){\mathcal{A}_d} [rrr] & TQ \ar[r]_{\mathcal{A}} &
\mathfrak{g} \ar[r]_{\exp} & G }\, ,\]
we can rely on explicit expressions for the exponential and logarithm, or we can rely on approximations to the exponential and logarithm.

The explicit formulas for the special Euclidean group are particularly useful for applying the theory of discrete connections to the geometric control of problems such as robotic manipulators, and clusters of satellites. In dealing with other configuration manifolds, approximants to the exponential and logarithm may be required due to the absence of explicit formulas.

Even when explicit formulas are available, it may be desirable to rely on a more computationally efficient approximation, such as the Cayley transformation, methods based on Pad\'e approximants (see, for example,~\cite{CaSi2001,Higham2001}), or Lie group techniques (see, for example,~\cite{CeIs2000,CeIs2001,ZaMu2001}). Clearly, these will yield a discrete connection that has an order of approximation equal to the lower of the two orders of approximation of the numerical schemes used for the exponential and the logarithm.

When used in the context of geometric control, high-order approximations to the continuous connection may not be necessary, since the optimal trajectory is often recomputed at each step, and in such instances, a low-order approximation suffices.

\subsection{Discrete Mechanical Connections and Discrete Lagrangians}\label{dcpb:subsec:discrete_mechanical_connction_discrete_lagrangian}\index{discrete connection!mechanical!from discrete Lagrangian}

We will introduce a discrete mechanical connection that is consistent with the structure of discrete variational mechanics, and will yield a discrete connection that has an order of approximation that is equation to that obtained from the discrete mechanics.

Consider a $G$-invariant $k$-th order discrete Lagrangian, $L_d:Q\times
Q\rightarrow\mathbb{R}$, which is to say that
\[ L_d(gq_0, gq_1) = L_d(q_0, q_1), \]
and
\[ L_d = L_d^{\operatorname{exact}} + \mathcal{O}(h^{k+1}),\]
where the exact discrete Lagrangian, $L_d^{\operatorname{exact}}:Q\times Q\rightarrow \mathbb{R}$, is given by
\[ L_d^{\operatorname{exact}}(q_0,q_1) = \int_0^h L(q_{01}(t),\dot q_{01}(t))dt.\]
Here, $q_{01}:[0,h]\rightarrow Q$ is the solution of the Euler--Lagrange equations with $q_{01}(0)=q_0$, and $q_{01}(h)=q_1$. The exact discrete Lagrangian is a generator of the symplectic flow, coming from the Jacobi solution of the Hamilton--Jacobi equation.

This $k$-th order discrete Lagrangian yields a $k$-th order accurate numerical update scheme, through the discrete Euler--Lagrange equations,
\[ D_2 L_d(q_0,q_1) + D_1 L_d(q_1, q_2) = 0,\]
which implicitly defines a discrete flow
$\Phi:(q_0,q_1)\mapsto(q_1,q_2)$. By pushing this numerical scheme
forward to $T^*Q$ using the discrete fiber derivative
$\mathbb{F}L_d:Q\times Q\rightarrow T^*Q$, which maps
$(q_0,q_1)\mapsto(q_0,-D_1 L_d(q_0,q_1))$, we can obtain a
Symplectic Partitioned Runge--Kutta scheme of the same order.

We also introduce the discrete momentum map, $J_d:Q\times
Q\rightarrow \mathfrak{g}^*$, given by
\[ \langle J_d(q_0,q_1),\xi \rangle = -D_1 L_d(q_0,q_1)\cdot
\xi_Q(q_0).\]

The discrete Lagrangian is $G$-invariant, which implies that for
any $\xi\in\mathfrak{g}$, we have,
\begin{align*}
L_d(q_0,q_1)&=L_d(\exp(\xi t)\cdot q_0, \exp(\xi t)\cdot q_1),\\
0&=\left.\frac{d}{dt}\right |_{t=0} L_d(\exp(\xi t)\cdot q_0,\exp
(\xi t)\cdot q_1)\\
&=D_1 L_d(q_0,q_1)\cdot \xi_Q(q_0) + D_2 L_d(q_0,q_1) \cdot
\xi_Q(q_1).
\end{align*}
If we restrict to the flow of the discrete Euler--Lagrange
equations, we have that
\[ (D_1 L_d(q_0,q_1)+D_2 L_dq(q_1,q_2))\cdot \xi_Q(q_1) = 0,\]
which upon substitution into the previous equation, yields
\begin{align*}
D_1 L_d(q_0,q_1)\cdot \xi_Q(q_0) - D_1 L_d(q_1,q_2)\cdot
\xi_Q(q_1) &= 0,\\
-D_1 L_d(q_1,q_2)\cdot\xi_Q(q_1) &= -D_1 L_d(q_0,q_1)
\cdot\xi_Q(q_0),\\
\langle J_d(q_1,q_2),\xi_Q(q_1) \rangle &= \langle
J_d(q_0,q_1),\xi_Q(q_0) \rangle,\\
\Phi^* J_d &= J_d.
\end{align*}
which is the statement of the discrete Noether theorem, that the
discrete momentum map is preserved by the discrete Euler--Lagrange
flow.

We note that the mechanical connection corresponds to a choice of
horizontal space corresponding to the zero momentum surface. That
is to say that the horizontal distribution corresponding to the
continuous mechanical connection is
\[ \operatorname{Hor}_q = \{v_q\in TQ\mid J(v_q)=0 \},
\]
and the discrete horizontal subspace corresponding to the discrete
mechanical connection is
\[ \operatorname{Hor}_q^d = \{ (q,q')\in Q\times Q \mid
J_d(q,q')=0 \}\]

\begin{remark}
For the discrete horizontal subspace we defined above to have the correct dimensionality,
the discrete Lagrangian needs to satisfy certain non-degeneracy conditions, which dictates the size of the neighborhood of the diagonal that the discrete connection is defined on.
\end{remark}

Since the continuous momentum map is preserved by the continuous
Euler--Lagrange flow, and the discrete momentum map is preserved
by the discrete Euler--Lagrange flow, it follows that the order of
approximation of the continuous mechanical connection by the
discrete mechanical connection is equal to the order of
approximation of the continuous Euler--Lagrange flow by the
discrete Euler--Lagrange flow. To construct a discrete mechanical
connection of a prescribed order, we use the following procedure.

\begin{enumerate}
\item Consider a $k$-th order $G$-invariant discrete Lagrangian,
$L_d:Q\times Q\rightarrow \mathbb{R}$,
\[ L_d = L_d^{\operatorname{exact}} + \mathcal{O}(h^{k+1}).\]
\item Construct the corresponding discrete momentum map,
$J_d(q_0,q_1)\rightarrow \mathfrak{g}^*$,
\[ \langle J_d(q_0,q_1),\xi\rangle = -D_1 L_d(q_0,q_1)\cdot
\xi_Q(q_0) .\]
\item Then, the $k$-th order discrete mechanical
connection is given implicitly by considering the condition for the discrete horizontal space,
\[ \mathcal{A}_d(q_0,q_1)=e\qquad\text{iff}\qquad J_d(q_0,q_1)=0,
\]
and then extending the construction to the domain of definition by $G$-equivariance.
\item More explicitly, given $(q_0,q_1)\in Q\times Q$, we consider
a local trivialization, in which $(q_0,q_1)=(x_0,g_0,x_1,g_1)$,
and we find the unique $g\in G$ such that
\[ J_d(x_0,g_0,x_1,g)=0.\]
Then, we have that
\[ \mathcal{A}_d(x_0,g_0,x_1,g) = e,\]
from which we conclude that
\begin{align*}
\mathcal{A}_d(x_0,g_0,x_1,g_1)
&= \mathcal{A}_d(x_0,g_0,x_1,g_1 g^{-1} g )\\
&= g_1 g^{-1}\cdot\mathcal{A}_d(x_0,g_0,x_1,g)\\
&= g_1 g^{-1}\, .
\end{align*}
\end{enumerate}

\section{Applications}\label{dcpb:sec:applications}
This section will sketch some of the applications of the
mathematical machinery of discrete connections and discrete
exterior calculus to problems in computational geometric
mechanics, geometric control theory, and discrete Riemannian
geometry.

\subsection{Discrete Lagrangian Reduction}\index{discrete connection!Lagrangian reduction}
Lagrangian reduction, which is the Lagrangian analogue of Poisson
reduction on the Hamiltonian side, is associated with the
reduction of Hamilton's variational principle for systems with
symmetry.

The variation of the action integral associated with a variation
in the curve can be expressed in terms of the Euler--Lagrange
operator, $\mathcal{EL}:T^{(2)}Q\rightarrow T^* Q$. When the
Lagrangian is $G$-invariant, the associated Euler--Lagrange
operator is $G$-equivariant, and this induces a reduced
Euler--Lagrange operator, $[\mathcal{EL}]_G:T^{(2)}Q/G\rightarrow
T^*Q/G$. The choice of a connection allows us to construct
intrinsic coordinates on $T^{(2)}/G$ and $T^*Q/G$, and the
representation of the reduced Euler--Lagrange operator in these
coordinates correspond to the Lagrange--Poincar\'e operator,
$\mathcal{LP}:T^{(2)}(Q/G)\times_{Q/G}
2\tilde{\mathfrak{g}}\rightarrow
T^*(Q/G)\oplus\tilde{\mathfrak{g}}^*$.

The reduced equations obtained by reduction tend to have
non-canonical symplectic structures. As such, na\"\i vely applying
standard symplectic algorithms to reduced equations can have
undesirable consequences for the longtime behavior of the
simulation, since it preserves the canonical symplectic form on
the reduced space, as opposed to the reduced (non-canonical)
symplectic form that is invariant under the reduced dynamics.

This sends a cautionary message, that it is important
to understand the reduction of discrete variational mechanics,
since applying standard numerical algorithms to the reduced
equations obtained from continuous reduction theory may yield undesirable results, inasmuch as long-term stability is concerned.

Discrete connections on principal bundles provide the appropriate
geometric structure to construct a discrete analogue of Lagrangian
reduction. We first introduce the discrete Euler--Lagrange
operator, which is constructed as follows.

\paragraph{Discrete Euler--Lagrange Operator.}\index{Euler--Lagrange!discrete operator} The discrete
Euler--Lagrange operator, $\mathcal{EL}_d:Q^3\rightarrow T^* Q$
satisfies the following property,
\[\mathbf{d}\mathfrak{S}_d(L_d)\cdot\delta\mathbf{q}=\sum
\mathcal{EL}_d(L_d)(q_{k-1},q_k,q_{k+1})\cdot\delta q_k.\] In
coordinates, the discrete Euler--Lagrange operator has the form
\[ \left[ D_2 L_d(q_{k-1},q_k) + D_1 L_d(q_k,q_{k+1}) \right
]dq_k.\]

\paragraph{Discrete Lagrange--Poincar\'{e} Operator.}\index{reduction!Lagrange--Poincar\'e}\index{Lagrange--Poincar\'e!discrete operator} The map
$\mathcal{EL}_d(L_d):Q^3\rightarrow T^*Q$, being $G$-equivariant,
induces a quotient map
\[[\mathcal{EL}_d(L_d)]_G:Q^3/G\rightarrow T^*Q/G,\]
which depends only on the reduced discrete Lagrangian
$l_d:(Q\times Q)/G\rightarrow\mathbb{R}.$ We can therefore
identify $[\mathcal{EL}_d(L_d)]_G$ with an operator
$\mathcal{EL}_d(l_d)$ which we call the {\bfi reduced discrete
Euler--Lagrange operator}\index{Euler--Lagrange!reduced!discrete}.

If in addition to the principal bundle structure, we have a
discrete principal connection as described in the previous
section, we can identify
\[Q^3/G\quad\text{with}\quad
(Q/G)^3\times_{Q/G} (\tilde{G}\oplus\tilde{G}).\]
The isomorphism between these two spaces is a consequence of Proposition~\ref{dcpb:prop:higher_order_isomorphism}, which is higher-order generalization of Proposition~\ref{dcpb:prop:isomorphism}. The discrete mechanical connection which was introduced in \S\ref{dcpb:subsec:discrete_mechanical_connction_discrete_lagrangian} is a particularly natural choice of discrete connection, since it does not require any \textit{ad hoc} choices, as it is constructed directly from the discrete Lagrangian.

Furthermore, each discrete $G$-valued connection $1$-form, $\mathcal{A}_d:Q\times Q\rightarrow G$, induces in the infinitesimal limit a continuous $\mathfrak{g}$-valued connection $1$-form, $\mathcal{A}:TQ\rightarrow
\mathfrak{g}$, as shown in \S\ref{s:infinitesimal_limit}. This continuous principal connection allows us to identify
\[T^*Q/G\quad\text{with}\quad
T^*(Q/G)\oplus\tilde{\mathfrak{g}}^*.\]

The discrete Lagrange--Poincar\'{e} operator,
$\mathcal{LP}_d(l_d):(Q/G)^3\times_{Q/G}
(\tilde{G}\oplus\tilde{G})\rightarrow T^* (Q/G) \oplus
\tilde{\mathfrak{g}}^*$, is obtained from the reduced discrete
Euler--Lagrange operator by making the identifications obtained
from the discrete connection structure.

The splitting of the range space of $\mathcal{LP}_d(l_d)$ as a
direct product (as in \S 3.3 of \cite{CeMaRa2001}) naturally induces a decomposition of the
discrete Lagrange--Poincar\'{e} operator,
\[\mathcal{LP}_d(l_d)=\operatorname{Hor}(\mathcal{LP}_d(l_d))\oplus
\operatorname{Ver}(\mathcal{LP}_d(l_d)),\] and this allows the
discrete reduced equations to be decomposed in horizontal and
vertical equations.

\subsection{Geometric Control Theory and Formations}\index{discrete connection!geometric control}
There are well-established control algorithms for actuating a
control system to achieve a desired reference configuration. In
many problems of practical interest, the actuation of the
mechanical system decomposes into shape and group variables in a
natural fashion.

A canonical example of this is a satellite in motion about the
Earth, where the orientation of the satellite is controlled by
internal rotors through the use of holonomy and geometric phases,
and the position is controlled by chemical propulsion.

In this example, the configuration space is $\SE(3)$, the group is
$\SO(3)$, and the shape space is $\mathbb{R}^3$. The group variable corresponds to the orientation, and the shape variable corresponds to the position. When given an
initial and desired configuration, it is desirable, while computing
the control inputs, to decompose the relative motion into a shape
component and a group component, so that they can be individually
actuated.

Since the discrete connection is used here to provide an efficient choice of local coordinates for optimal control, the discrete connection is most naturally obtained by exponentiating the continuous connection, in the manner described in \S\ref{dcpb:subsec:exponentiated_discrete_connections}, in conjunction with the explicit formulas for the exponential and logarithm for $\SE(3)$. The natural choice of the continuous connection is one in which the horizontal space is given by the momentum surface corresponding to the current value of the momentum.

\begin{figure}[htbp]
\begin{center}
\begin{minipage}{0.4\textwidth}
\begin{center}
\includegraphics[scale=0.5]{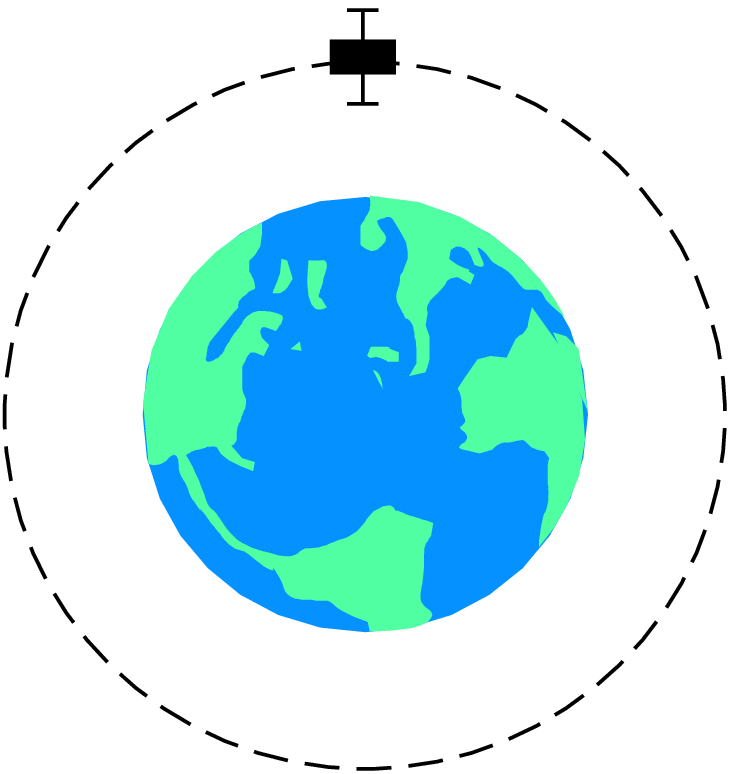}\\
Initial configuration
\end{center}
\end{minipage}\qquad
\begin{minipage}{0.4\textwidth}
\begin{center}
\includegraphics[scale=0.5]{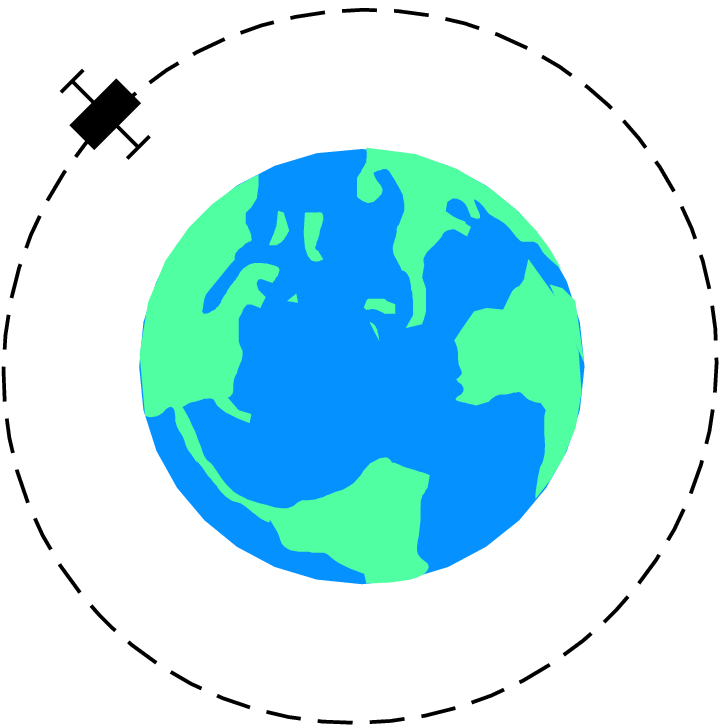}\\
Desired configuration
\end{center}
\end{minipage}
\end{center}
\caption{\label{dcpb:fig:control}Application of discrete connections to control.}
\end{figure}

To illustrate why it may not be desirable from a control-theoretic
point of view to decompose the space using a trivial connection,
consider a satellite that is in a tidally locked orbit about the
Earth, with the initial and desired configuration as illustrated in Figure~\ref{dcpb:fig:control}.

Here, if we choose a trivial connection, then the relative group
element would be a rotation by $\pi/4$, but this choice is
undesirable, since the motion is tidally locked, and moving the
center of mass to the new location would result in a shift in the
orientation by precisely the desired amount. In this case, the
optimal control input should therefore only actuate in the shape
variables, and the relative group element assigned to this pair of
configurations should be the identity element.

As such, the extension of mechanically relevant connections to
pairs of points in the configuration space with finite separation,
through the use of a discrete connection, can be of immense value in
geometric control theory.

Similarly, in the case of formations, discrete connections allow
for the orientation coordination problem to be handled in a more
efficient manner, by taking into account the dynamic coupling of
the shape and group motions automatically through the use of the discrete
mechanical connection.

\subsection{Discrete Levi-Civita Connection}\label{dcpb:subsec:levi_civita}\index{discrete connection!Levi-Civita}

Vector bundle connections can be cast in the language of
connections on principal bundles by considering the frame bundle
consisting of oriented orthonormal frames over the manifold $M$,
which is a principal $\SO(n)$ bundle, as originally proposed by
\cite{Ca1983,Ca2001}. For related work on discrete connections on triangulated manifolds with applications to algebraic topology and the computation of Chern classes, please see \cite{Novikov2003}.

To construct our model of a discrete Riemannian manifold, we first
trivialize the frame bundle to yield $\SO(n)\times M$. Then,
$Q=\SO(n)\times M$, and $G=\SO(n)$.

Here, we introduce the notion of a semidiscretized principal
bundle, where the shape space, $S=Q/G$, is discretized as a
simplicial complex $K$, and the structure group $G$ remains
continuous. In this context, the semidiscretization of the
trivialization of the frame bundle is given by $Q=\SO(n)\times K$.

A discrete connection is a map $\mathcal{A}_d:Q\times Q\rightarrow G$,
and we can construct a candidate for the Levi-Civita connection on
a simplicial complex $K$, using the discrete analogue of the frame
bundle described above. However, we will first introduce the notion of a discrete Riemannian manifold.

\begin{definition}
A \textbf{discrete dual Riemannian manifold}\index{discrete Riemannian manifold} is a simplicial
complex where each $n$-simplex $\sigma^n$ is endowed with a
constant Riemannian metric tensor $g$, such that the restriction
of the metric tensor to a common face with an adjacent $n$-simplex
is consistent.
\end{definition}

This is referred to as a discrete dual Riemannian manifold as we
can equivalently think of associating a Riemannian metric tensor
to each dual vertex, and as we shall see, by adopting Cartan's
method of orthogonal frames (see, for example, \cite{Ca1983, Ca2001}), the connection is a
$\SO(n)$-valued discrete dual $1$-form, and the curvature is a
$\SO(n)$-valued discrete dual $2$-form.

For each $n$-simplex $\sigma^n$, consider an invertible
transformation $f$ of $\RR^n$ such that $f^*g=I.$ In the
orthonormal space, we have a normal operator that maps a $(n-1)$-dimensional subspace to a generator of the orthogonal complement,
denoted by $\perp$. Then, we obtain a normal operator on the faces
of $\sigma^n$ by making the following diagram commute.

\[\xymatrix@!0@R=3cm@C=4cm{
\includegraphics[scale=0.2]{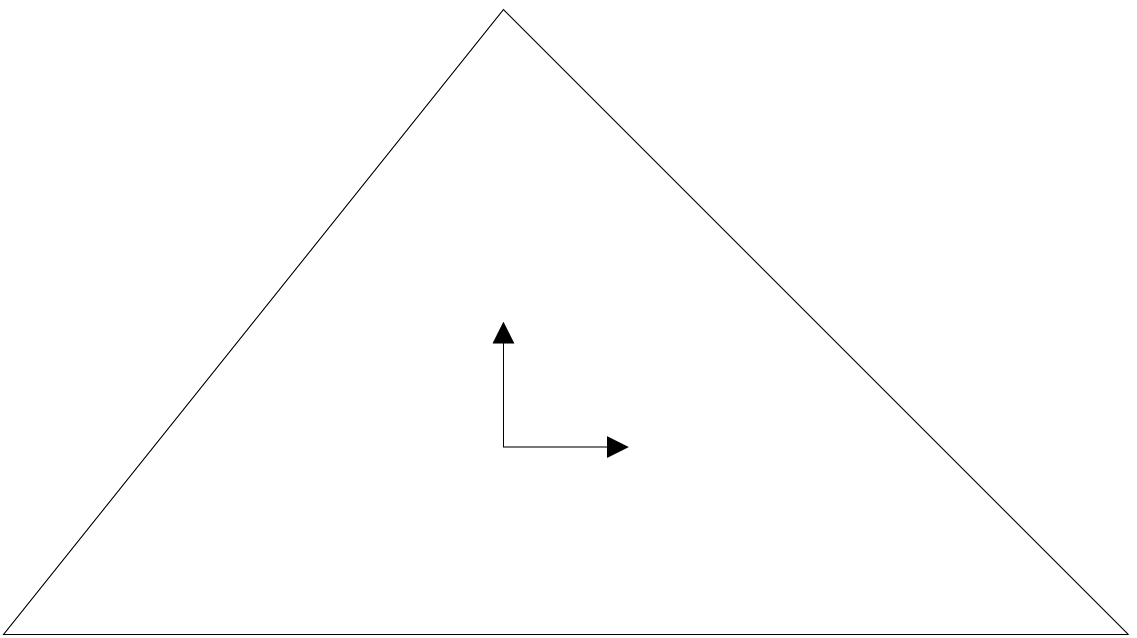} \ar[r]^f \ar[d]_{\perp} &
\includegraphics[scale=0.15]{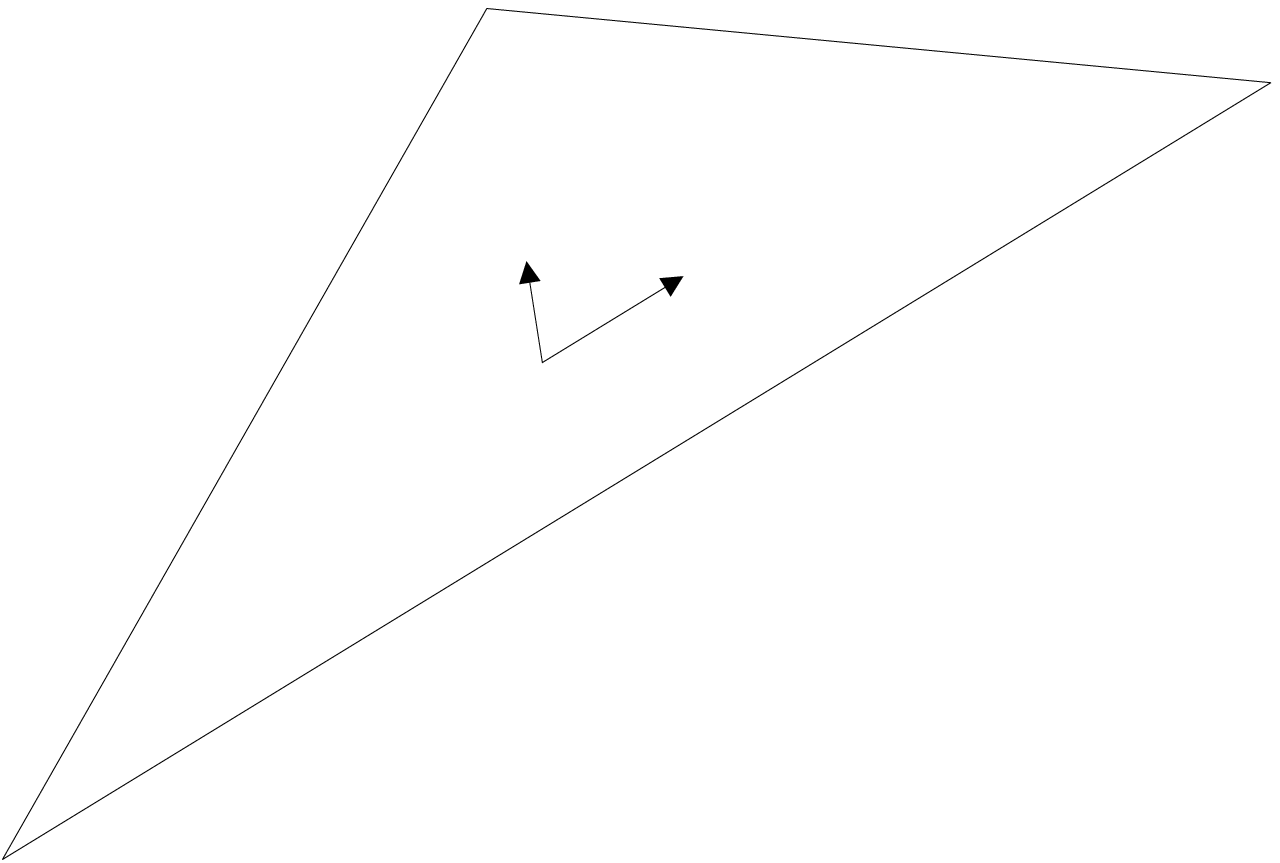} \ar[d]^{\perp}\\
\includegraphics[scale=0.2]{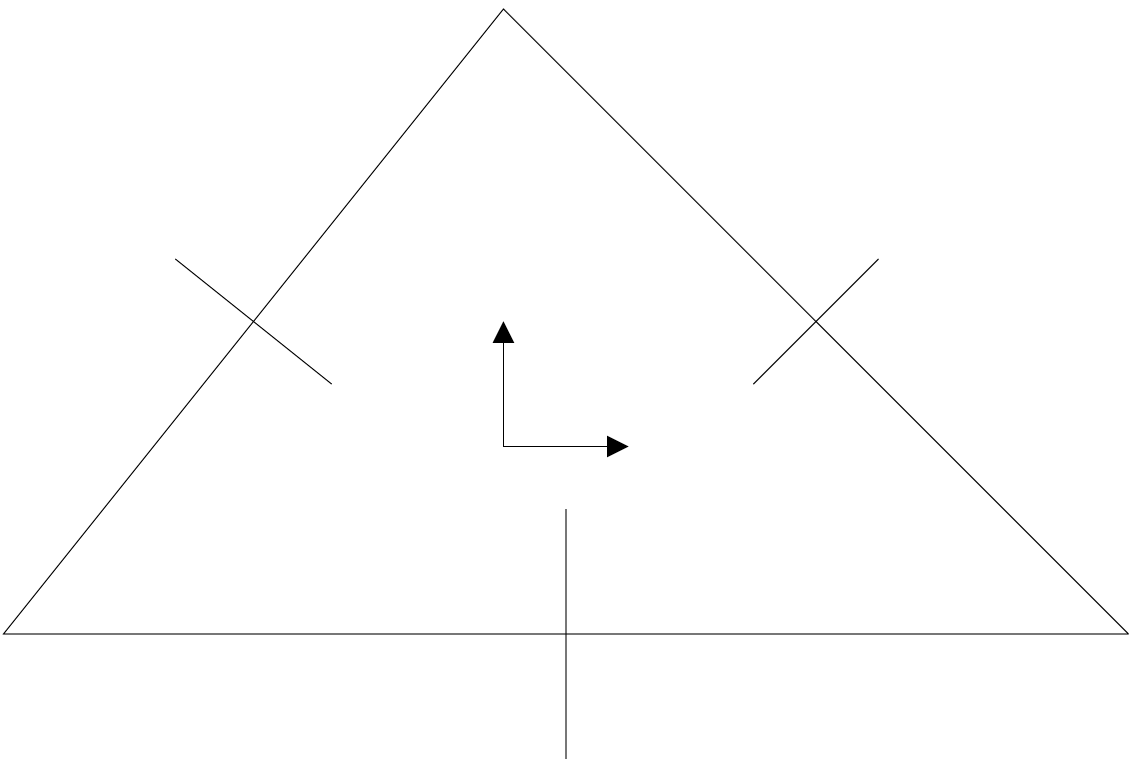} \ar[r]^f &
\includegraphics[scale=0.15]{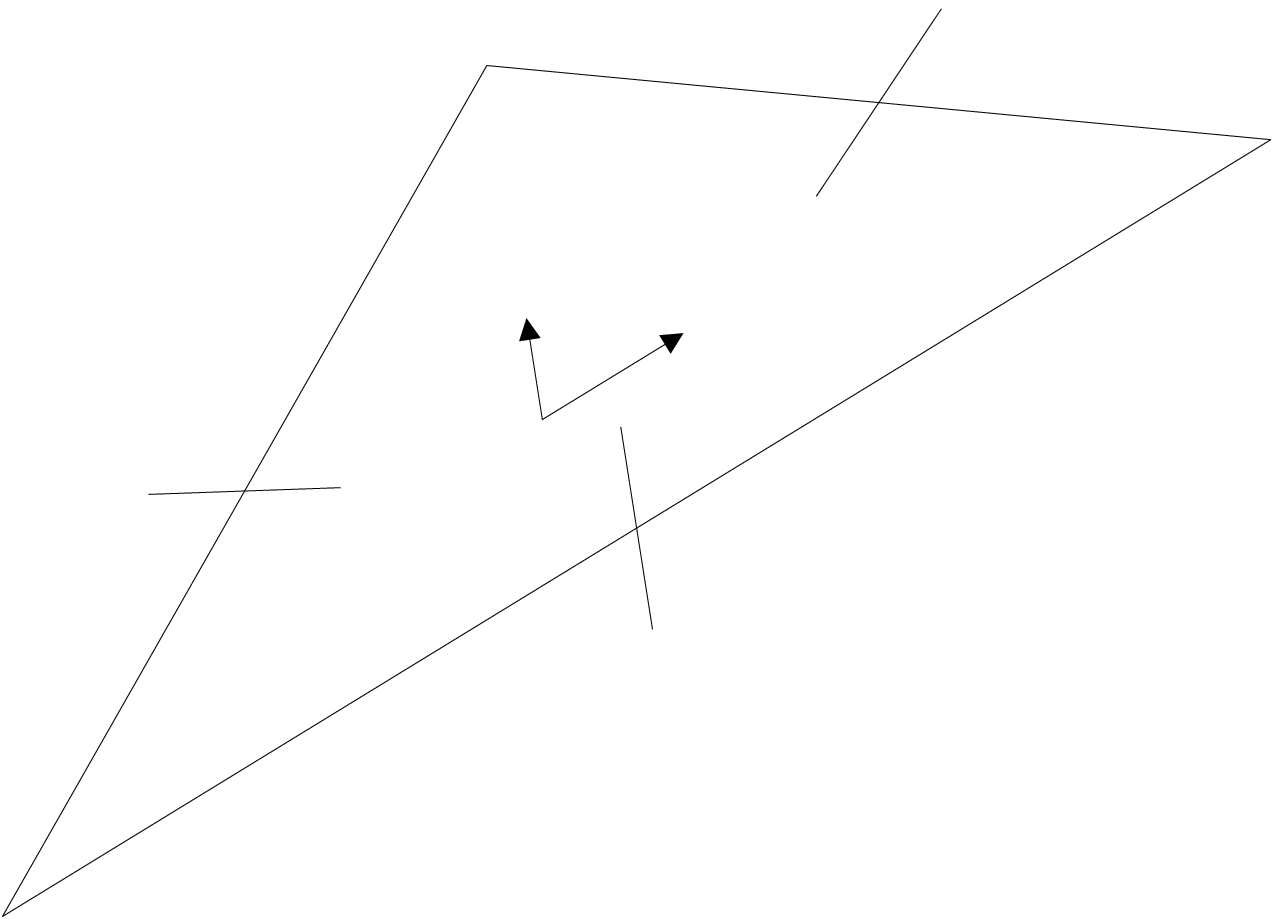}}\]

The coordinate axes in the diagram represent the normalized
eigenvectors of the metric, scaled by their respective
eigenvalues.

The local representation of the discrete connection is given by
\[\mathcal{A}_d((\sigma^n_0, R_0),(\sigma^n_1, R_1))=R_1
A(\sigma^n_0,\sigma^n_1)R_0^{-1},\] and so the discrete connection
is uniquely defined if we specify $A(\sigma^n_0,\sigma^n_1)$,
where $\sigma^n_0$ and $\sigma^n_1$ are adjacent $n$-simplices.
Since they are adjacent, they share a $(n-1)$-simplex, denoted
$\sigma^{n-1}$. In particular, this can then be thought of as a
$\SO(n)$-valued discrete dual $1$-form, since to each dual
$1$-cell, $\star\sigma^{n-1}$, we associate an element of $\SO(n)$.

This element of $\SO(n)$ is computed as follows.
\begin{enumerate}
\item In each of the $n$-simplices, we have a normal direction
associated with $\sigma^{n-1}$, denoted by
\[\perp(\sigma^{n-1},\sigma^n_i)\in\RR^n.\]
\item If these two normal directions are parallel, we set
\[\langle A,\star\sigma^{n-1}\rangle =I,\]
otherwise, we continue. \item Construct the $(n-2)$-dimensional
hyperplane $P^{n-2}$, given by the orthogonal complement to the
span of the two normal directions.
\[P^{n-2}=\perp(\operatorname{span}(\perp(\sigma^{n-1},\sigma^n_0),
\perp(\sigma^{n-1},\sigma^n_1))).\] \item If $\star\sigma^{n-1}$
is oriented from $\sigma^n_0$ to $\sigma^n_1$, we set
\[\langle A,\star\sigma^{n-1}\rangle = \{R\in \SO(n) \mid
R|_{P^{n-2}}= I_{P^{n-2}}, R
(\perp(\sigma^{n-1},\sigma^n_0))=\perp(\sigma^{n-1},\sigma^n_1)
\}.\]
\end{enumerate}

The curvature of this discrete Levi-Civita connection\index{discrete connection!Levi-Civita!curvature} is then a $\SO(n)$-valued discrete dual $2$-form. There is however the
curious property that the boundary operator for a dual cell complex may not necessarily agree with the standard notion of boundary, since that may not be expressible in terms of a chain in the dual cell complex. This is primarily an issue on the boundary of the simplicial complex, and if we are in the interior, this is not a problem.

\begin{figure}[htbp]
\singlespace
\begin{center}
\begin{minipage}{0.19\textwidth}
\begin{center}
\includegraphics[scale=0.4,clip=true]{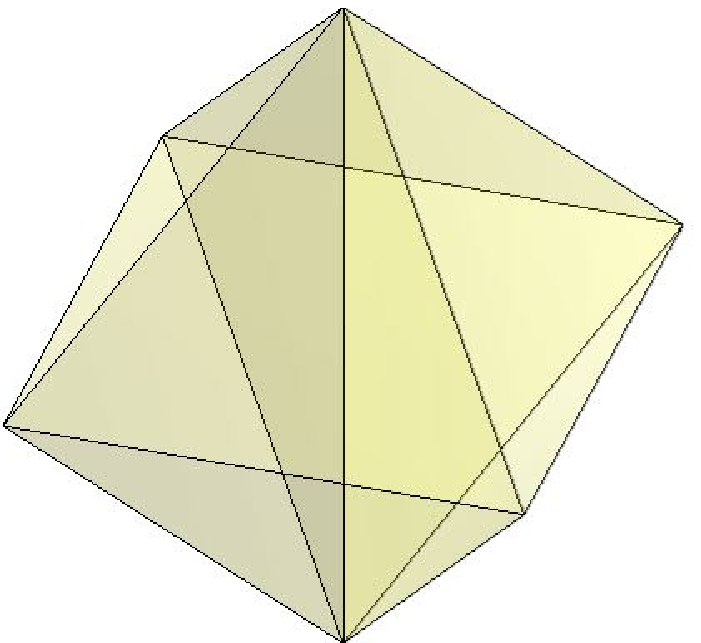}\\
Simplicial Complex,\\
$K$
\end{center}
\end{minipage}
\begin{minipage}{0.19\textwidth}
\begin{center}
\includegraphics[scale=0.4,clip=true]{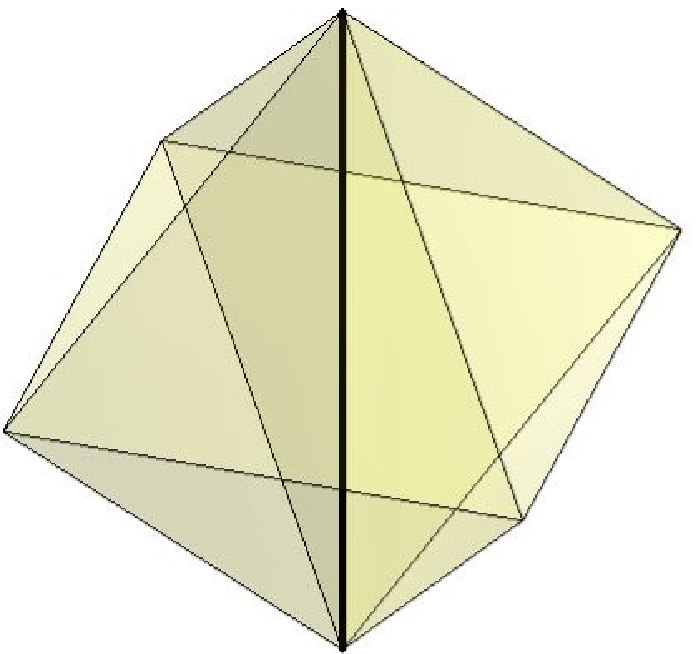}\\
primal $(n-2)$-simplex,\\
$\sigma^{n-2}$
\end{center}
\end{minipage}
\begin{minipage}{0.19\textwidth}
\begin{center}
\includegraphics[scale=0.4,clip=true]{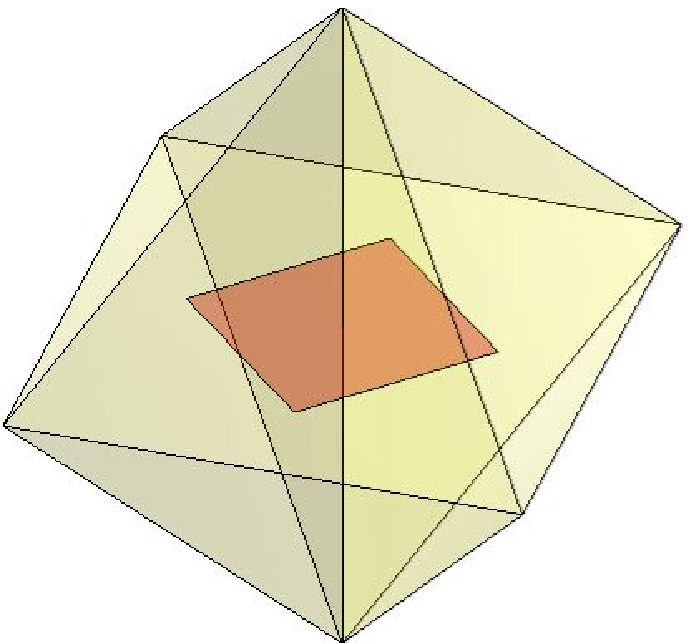}\\
dual $2$-cell,\\
$\star\sigma^{n-2}$
\end{center}
\end{minipage}
\begin{minipage}{0.19\textwidth}
\begin{center}
\includegraphics[scale=0.4,clip=true]{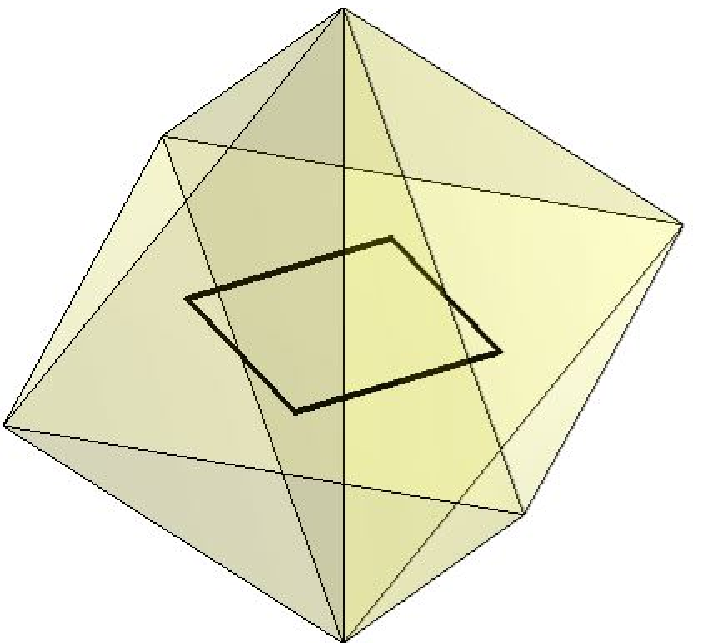}\\
dual $1$-chain,\\
$\partial\star\sigma^{n-2}$
\end{center}
\end{minipage}
\begin{minipage}{0.19\textwidth}
\begin{center}
\includegraphics[scale=0.4,clip=true]{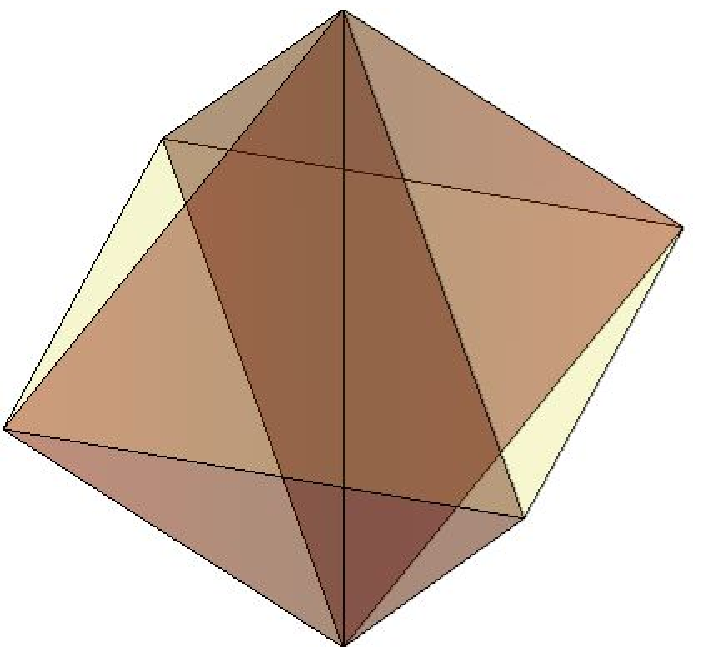}\\
primal $(n-1)$-chain,\\
$\star\partial\star\sigma^{n-2}$
\end{center}
\end{minipage}
\end{center}
\caption{\label{dcpb:fig:curvature_2_form}Discrete curvature as a discrete dual $2$-form.}
\end{figure}

Since curvature is a dual $2$-form, it is associated with the dual
of a codimension-two simplex, given by $\star\sigma^{n-2}$. Consider
the example illustrated in Figure~\ref{dcpb:fig:curvature_2_form}.

The curvature\index{curvature!Levi-Civita connection|see{connection, Levi-Civita, curvature}} $\mathcal{B}$ of the discrete Levi-Civita connection
is given by
\[\langle\mathcal{B},\star\sigma^{n-2}\rangle = \langle\d A, \star\sigma^{n-2} \rangle = \langle A,
\partial\star\sigma^{n-2} \rangle.\]
As can be seen from the geometric region
$\star\partial\star\sigma^{n-2}$, the curvature associated with
$\star\sigma^{n-2}$ is given by the ordered product of the
connection associated with the dual cells, $\star\sigma^{n-1}$,
where $\sigma^{n-1}\succ\sigma^{n-2}$. This also suggests that we
can also think of the discrete connection as a primal $(n-1)$-form,
and the curvature as a primal $(n-2)$-form, where the curvature is
obtained from the connection using the codifferential.

When the group $G$ is nonabelian, we see that the curvature is
only defined up to conjugation, since we need to specify a dual
vertex on the dual one-chain $\partial\star\sigma^{n-2}$ from which
to start composing group elements. To make this well-defined, we
can adopt the approach used in defining the simplicial cup
product, and assume that there is a partial ordering on dual
vertices, which would make the curvature unambiguous.

As a quality measure for simplicial triangulations of a Riemannian
manifold, having the curvature defined up to conjugation may be
sufficient if we have a norm on $\SO(n)$ which is invariant under
conjugation. As an example, taking the logarithm to get an element
of the Lie algebra $\mathfrak{so}(n)$, and then using a norm on this vector
space yields a conjugation-invariant norm on $\SO(n)$. This allows
us to detect regions of the mesh with high curvature, and
selectively subdivide the triangulation in such regions.

Similarly, we can define a discrete primal Riemannian manifold,
where the Riemannian metric tensor is associated with primal
vertices, and the connection is a $G$-valued primal $1$-form, and
the curvature is a $G$-valued primal $2$-form.

\paragraph{Abstract Simplicial Complex with a Local Metric.}\index{metric!local}
In \cite{DeHiLeMa2003}, the notion of an abstract simplicial complex with a
local metric, defined on pairs of vertices that are adjacent, was introduced.

In this situation, we can compute the curvature around a loop in the
mesh using local embeddings. We start with an initial $n$-simplex,
which we endow with an orthonormal frame. By locally embedding
adjacent $n$-simplices into Euclidean space, and parallel
transporting the orthonormal frame, we will eventually transport
the frame back to the initial simplex.

The relative orientation between the original frame and the
transported frame yields the integral of the curvature of the
surface which is bounded by the traversed curve. This results from
a simple application of the Generalized Stokes' theorem, and the
fact that the curvature is given by the exterior derivative of the
connection $1$-form.

\section{Conclusions and Future Work}
We have introduced a complete characterization of discrete connections, in terms of horizontal and vertical spaces, discrete connection $1$-forms, horizontal lifts, and splittings of the discrete Atiyah sequence. Geometric structures that can be derived from a given discrete connection have been discussed, including continuous connections, an extended pair groupoid composition, and higher-order analogues of the discrete connection. In addition, we have explored computational issues, such as order of accuracy, and the construction of discrete connections from continuous connections.

Applications to discrete reduction theory, geometric control theory, and discrete geometry, have also been discussed, and it would be desirable to systematically apply the machinery of discrete connections to these problems.

In addition, connections play a crucial role in representing the nonholonomic constraint distribution in nonholonomic mechanics, particularly when considering nonholonomic mechanical systems with symmetry, wherein the nonholonomic connection enters (see, for example, \cite{Bloch2003}).
There has been recent progress  on constructing nonholonomic integrators in the work of \cite{Cortes2002} and \cite{McLaPe2003}, but an intrinsically discrete notion of a connection remains absent from their work, and they do not consider the role of symmetry reduction in discrete nonholonomic mechanics.

It would be very interesting to apply the general theory of discrete connections on principal bundles to nonholonomic mechanical systems with symmetry, and to cast the notion of a discrete nonholonomic constraint distribution and the nonholonomic connection in the language of discrete connections, and thereby develop a discrete theory of nonholonomic mechanics with symmetry. This would be particularly important for the numerical implementation of geometric control algorithms.

The role of discrete connections in the study of discrete geometric phases would also be an area worth pursuing. A discrete analogue of the rigid-body phase formula, that involves the discrete mechanical connection, that is exact for rigid-body simulations that use discrete variational mechanics, would yield significant insights into the geometric structure-preservation properties of variational integrators. In particular, it would provide much needed insight into how discretization interacts with geometric phases, and yield an understanding how much of the phase drift observed in a numerical simulation is due to the underlying geometry of the mechanical system, and how much is due to the process of discretizing the system.

\bibliographystyle{plainnat}
\bibliography{umich_dcpb}

\begin{thebibliography}{30}
\providecommand{\natexlab}[1]{#1}
\providecommand{\url}[1]{\texttt{#1}}
\expandafter\ifx\csname urlstyle\endcsname\relax
  \providecommand{\doi}[1]{doi: #1}\else
  \providecommand{\doi}{doi: \begingroup \urlstyle{rm}\Url}\fi

\bibitem[Almeida and Molino(1985)]{AlMo1985}
R.~Almeida and P.~Molino.
\newblock Suites d'{A}tiyah et feuilletages transversalement complets.
\newblock \emph{C. R. Acad. Sci. Paris S\'er. I Math.}, 300\penalty0
  (1):\penalty0 13--15, 1985.

\bibitem[Atiyah(1957)]{At1957}
{M. F.} Atiyah.
\newblock Complex analytic connections in fibre bundles.
\newblock \emph{Trans. Amer. Math. Soc.}, 85:\penalty0 181--207, 1957.

\bibitem[Berry(1990)]{Berry1990}
M.~Berry.
\newblock Anticipations of the geometric phase.
\newblock \emph{Phys. Today}, pages 34--40, December 1990.

\bibitem[Bloch(2003)]{Bloch2003}
{A. M.} Bloch.
\newblock \emph{Nonholonomic Mechanics and Control}.
\newblock Interdisciplinary Applied Mathematics. Springer-Verlag, 2003.

\bibitem[Cardoso and {Silva Leite}(2001)]{CaSi2001}
{J. R.} Cardoso and F.~{Silva Leite}.
\newblock Theoretical and numerical considerations about {P}ad\'e approximants
  for the matrix logarithm.
\newblock \emph{Linear Algebra Appl.}, 330\penalty0 (1-3):\penalty0 31--42,
  2001.

\bibitem[Cartan(1983)]{Ca1983}
\'E. Cartan.
\newblock \emph{Geometry of Riemannian spaces}.
\newblock Math. Sci. Press, Brookline, CA, 1983.
\newblock (translated from French).

\bibitem[Cartan(2001)]{Ca2001}
\'E. Cartan.
\newblock \emph{Riemannian Geometry in an Orthogonal Frame}.
\newblock World Scientific, 2001.
\newblock (translated from French).

\bibitem[Celledoni and Iserles(2000)]{CeIs2000}
E.~Celledoni and A.~Iserles.
\newblock Approximating the exponential from a {L}ie algebra to a {L}ie group.
\newblock \emph{Math. Comp.}, 69\penalty0 (232):\penalty0 1457--1480, 2000.

\bibitem[Celledoni and Iserles(2001)]{CeIs2001}
E.~Celledoni and A.~Iserles.
\newblock Methods for the approximation of the matrix exponential in a
  lie-algebraic setting.
\newblock \emph{IMA J. Num. Anal.}, 21\penalty0 (2):\penalty0 463--488, 2001.

\bibitem[Cendra et~al.(2001)Cendra, Marsden, and Ratiu]{CeMaRa2001}
H.~Cendra, {J. E.} Marsden, and {T. S.} Ratiu.
\newblock Lagrangian reduction by stages.
\newblock \emph{Mem. Amer. Math. Soc.}, 152\penalty0 (722), 2001.

\bibitem[Cort\'es(2002)]{Cortes2002}
J.~Cort\'es.
\newblock \emph{Geometric, Control and Numerical Aspects of Nonholonomic
  Systems}, volume 1793 of \emph{Lecture Notes in Mathematics}.
\newblock Springer-Verlag, 2002.

\bibitem[Desbrun et~al.(2003)Desbrun, Hirani, Leok, and Marsden]{DeHiLeMa2003}
M.~Desbrun, {A. N.} Hirani, M.~Leok, and {J. E.} Marsden.
\newblock Discrete exterior calculus.
\newblock (in preparation), 2003.

\bibitem[Goldreich and Toomre(1969)]{GoTo1969}
P.~Goldreich and A.~Toomre.
\newblock Some remarks on polar wandering.
\newblock \emph{J. Geophys. Res.}, 10:\penalty0 2555--2567, 1969.

\bibitem[Higham(2001)]{Higham2001}
{N. J.} Higham.
\newblock Evaluating {P}ad\'e approximants of the matrix logarithm.
\newblock \emph{SIAM J. Matrix Anal. Appl.}, 22\penalty0 (4):\penalty0
  1126--1135 (electronic), 2001.

\bibitem[Kobayashi and Nomizu(1963)]{KoNo1963}
S.~Kobayashi and K.~Nomizu.
\newblock \emph{Foundations of Differential Geometry}, volume~1.
\newblock Wiley, 1963.

\bibitem[Leok(1998)]{Le1998}
M.~Leok.
\newblock A mathematical model of true polar wander.
\newblock Caltech SURF Report, 1998.

\bibitem[Mackenzie(1995)]{Mackenzie1995}
{K. C. H.} Mackenzie.
\newblock Lie algebroids and {L}ie pseudoalgebras.
\newblock \emph{Bull. London Math. Soc.}, 27\penalty0 (2):\penalty0 97--147,
  1995.

\bibitem[Marsden(1994)]{Marsden1994}
{J. E.} Marsden.
\newblock Geometric mechanics, stability, and control.
\newblock In L.~Sirovich, editor, \emph{Applied Mathematical Sciences}, volume
  100, pages 265--291. Springer-Verlag, 1994.

\bibitem[Marsden(1997)]{Marsden1997}
{J. E.} Marsden.
\newblock Geometric foundations of motion and control.
\newblock In \emph{Motion, Control and Geometry}, pages 3--19. National Academy
  Press, 1997.

\bibitem[Marsden and Ratiu(1999)]{MaRa1999}
{J. E.} Marsden and {T. S.} Ratiu.
\newblock \emph{Introduction to Mechanics and Symmetry}, volume~17 of
  \emph{Texts in Applied Mathematics}.
\newblock Springer-Verlag, second edition, 1999.

\bibitem[Marsden et~al.(1990)Marsden, Montgomery, and Ratiu]{MaMoRa1990}
{J. E.} Marsden, R.~Montgomery, and {T. S.} Ratiu.
\newblock Reduction, symmetry, and phases in mechanics.
\newblock \emph{Mem. Amer. Math. Soc.}, 88\penalty0 (436):\penalty0 1--110,
  1990.

\bibitem[McLachlan and Perlmutter(2003)]{McLaPe2003}
R.~McLachlan and M.~Perlmutter.
\newblock Integrators for nonholonomic mechanical systems.
\newblock (preprint), 2003.

\bibitem[Montgomery(1991)]{Montgomery1991}
R.~Montgomery.
\newblock How much does a rigid body rotate? {A} {B}erry's phase from the
  eighteenth century.
\newblock \emph{Amer. J. Phys.}, 59:\penalty0 394--398, 1991.

\bibitem[Novikov(2003)]{Novikov2003}
{S. P.} Novikov.
\newblock Discrete connections on the triangulated manifolds and difference
  linear equations.
\newblock arXiv, math-ph/0303035, 2003.

\bibitem[Olver(2001)]{Olver2001}
{P. J.} Olver.
\newblock Geometric foundations of numerical algorithms and symmetry.
\newblock \emph{Appl. Algebra Engrg. Comm. Comput.}, 11\penalty0 (5):\penalty0
  417--436, 2001.
\newblock Special issue ``Computational geometry for differential equations''.

\bibitem[Shapere and Wilczek(1989)]{ShWi1989}
A.~Shapere and F.~Wilczek.
\newblock \emph{Geometric Phases in Physics}, volume~5 of \emph{Advanced Series
  in Mathematical Physics}.
\newblock World Scientific, 1989.

\bibitem[Steenrod(1951)]{Steenrod1951}
N.~Steenrod.
\newblock \emph{The Topology of Fibre Bundles}.
\newblock Princeton University Press, 1951.

\bibitem[Vedral(2003)]{Vedral2003}
V.~Vedral.
\newblock Geometric phases and topological quantum computation.
\newblock \emph{Int. J. Quantum Inform.}, 1\penalty0 (1):\penalty0 1--23, 2003.

\bibitem[Walsh and Sastry(1993)]{WaSa1993}
G.~Walsh and S.~Sastry.
\newblock On reorienting linked rigid bodies using internal motions.
\newblock \emph{IEEE Trans. Robotic Autom.}, 11:\penalty0 139--146, 1993.

\bibitem[Zanna and Munthe-Kaas(2001/02)]{ZaMu2001}
A.~Zanna and {H. Z.} Munthe-Kaas.
\newblock Generalized polar decompositions for the approximation of the matrix
  exponential.
\newblock \emph{SIAM J. Matrix Anal. Appl.}, 23\penalty0 (3):\penalty0 840--862
  (electronic), 2001/02.

\end{thebibliography}

\end{document}